\documentclass[11pt,
]{amsart}
\usepackage{
amssymb, graphicx
}
\usepackage{mathrsfs} 
\usepackage{amsthm}
\usepackage{graphicx,color}
\usepackage[dvipsnames]{xcolor}
\usepackage{amsmath}
\allowdisplaybreaks[4]
\newtheorem{theorem}{Theorem}
\newtheorem{lemma}{Lemma}

\newtheorem{definition}{Definition}

\theoremstyle{remark}

\date{\today}

\title[An IBVP for the Maxwell's equations]{An inverse boundary value problem for the Maxwell's equations with partial data}

  \author[J. Zhai]{Jian Zhai}
\address{School of Mathematical Sciences,
  Fudan University, Shanghai 200433, China
  (\tt{jianzhai@fudan.edu.cn}).}
    \thanks{J. Zhai is supported by National Key Research and Development Programs of China (No. 2023YFA1009103), NSFC(No. 12471396), Science and Technology Commission of Shanghai Municipality (23JC1400501)}

\begin{document}
\begin{abstract}
We consider an inverse boundary problem for the dynamical Maxwell's equations. We show that the electric permittivity, conductivity, and magnetic permeability can be uniquely determined locally if there is a strictly convex foliation with respect to the wave speed.
\end{abstract}
\maketitle

\section{Introduction}

Consider the Maxwell's equations
\begin{equation}\label{maineq}
\mu\partial_t\mathcal{H}=-\nabla\times \mathcal{E},\quad \varepsilon\partial_t\mathcal{E}+\sigma \mathcal{E}=\nabla\times \mathcal{H}\quad \text{in }(0,T)\times\Omega,
\end{equation}
describing the traveling of electromagnetic waves in a bounded smooth domain $\Omega\subset\mathbb{R}^3$. Here $\mathcal{E}(t,x)$ is the electric field and $\mathcal{H}(t,x)$ is the magnetic field.
We assume that the electric permittivity $\varepsilon>0$, magnetic permeability $\mu>0$ and electric conductivity $\sigma\geq 0$ are all smooth functions in $\overline{\Omega}$.

The initial boundary value problem for the Maxwell's equations would be the equation \eqref{maineq} supplemented by the initial conditions
\begin{equation}\label{initialcondition}
\mathcal{E}(0,x)=0,\quad \mathcal{H}(0,x)=0
\end{equation}
and the boundary condition
\begin{equation}\label{boundarycondition}
\nu\times \mathcal{E}\vert_{(0,T)\times\partial\Omega}=f,
\end{equation}
where $\nu$ is the unit outer normal to $\partial\Omega$.

The inverse problem associated with the initial boundary value problem \eqref{maineq} \eqref{initialcondition} and \eqref{boundarycondition} is to determine $\varepsilon$, $\mu$ and $\sigma$ from the impedance map
\[
\Lambda_{\varepsilon,\mu,\sigma}: f\mapsto \nu\times\mathcal{H}\vert_{(0,T)\times\partial\Omega}.
\]

Tnverse boundary problems for electrodynamics have been studied since 1930s \cite{langer1933inverse,slichter1933inverse}, where the one-dimensional case was considered. For higher dimensional problems in frequency domain, we refer to \cite{ola1993inverse,ola1996electromagnetic,kenig2011inverse,caro2011inverse,brown2016uniqueness}. For problems in time domain, most results considered the $\sigma=0$ case, and the \textit{Boundary Control Method} (\cite{belishev1987approach,belishev2007recent}) can be applied.  The uniqueness of $\varepsilon\mu$ was established in \cite{belishev2001dynamical,belishev2000reconstruction}. 
For simultaneous recovery of $\varepsilon$ and $\mu$, we refer to \cite{belishev2004uniqueness,kurylev2006maxwell,demchenko2012dynamical}. In \cite{belishev2004uniqueness}, the problem is reduced to the one in frequency domain so the results in \cite{ola1993inverse,kenig2011inverse} can be applied. The work \cite{belishev2000reconstruction,demchenko2012dynamical} can even deal with the partial data problem in a time-optimal way. The $\sigma\neq 0$ case was studied in \cite{romanov2002inverse}, where the problem of recovering $\varepsilon,\mu,\sigma$ was reduced to certain geometric inverse problems. We will essentially adopt the same strategy here but taking a more microlocal perspective close in spirit to \cite{stefanov2018inverse}, and use the fact that all the associated geometric inverse problems have local uniqueness properties. More specifically speaking, the paper \cite{romanov2002inverse} proved the uniqueness of $\varepsilon\mu$ and $\frac{\sigma}{\varepsilon}$ from full boundary measurements under the assumption that $(\Omega,\varepsilon\mu\mathrm{d}s^2)$, considered to be a Riemannian manifold with boundary, is simple. Here  $\mathrm{d}s^2$ is the Euclidean metric. Recall that a Riemannian manifold $(\Omega,g)$ with boundary is called simple if any two points $x,y\in \Omega$ can be joined by a unique geodesic. However, the geometric inverse problem associated to the third parameter is more complicated, and a general uniqueness result is left open in \cite{romanov2002inverse}. We will see that this third geometric inverse problem is the so-called \textit{transverse ray transform} \cite{Shara}. The author of this article and collaborator proved the injectivity of the local transverse ray transform in \cite{uhlmann2024invertibility1,uhlmann2024invertibility2}, which can be applied here.\\
 
We consider the partial data problem in this article. For this, we first review the notion of a convex foliation with respect to the metric $g=c^{-2}\mathrm{d}s^2$ (cf. \cite{stefanov2017local}). Here $c=\frac{1}{\sqrt{\varepsilon\mu}}$ is the speed of the electromagnetic waves. Let $\widetilde{\Omega}$ be an open domain containing $\overline{\Omega}$, and extend the coefficients there in a smooth way. 
\begin{definition}
Let $\kappa:\widetilde{\Omega}\rightarrow\mathbb{R}$ be a smooth function for which the level sets $\kappa^{-1}(q),q\leq 1$, restricted to $\overline{\Omega}$, are strictly convex viewed from $\kappa^{-1}((-\infty,q))$ w.r.t. $g$; $\mathrm{d}\kappa\neq 0$ on these level sets, and $\kappa^{-1}(0)\cap\overline{\Omega}\subset\partial\Omega$. We call $\kappa^{-1}(q)\cap\overline{\Omega}$, $q\in [0,1]$, a strictly convex foliation of $\kappa^{-1}(q)\cap\overline{\Omega}$ with respect to $g$.
\end{definition}

Recall that an oriented hypersurface $S$ is strictly convex w.r.t. $g$, if the second fundamental form on $S$ is positive.  Now let $\Gamma=\kappa^{-1}([0,1))\cap\partial\Omega$, and define the local impedance map
\[
\Lambda_{\varepsilon,\mu,\sigma}^\Gamma:f\in C_c^\infty((0,T)\times \Gamma)\mapsto \nu\times\mathcal{H}\vert_{(0,T)\times \Gamma}
\]
The partial data problem is the recover $\varepsilon,\mu,\sigma$ in a neighborhood of $\Gamma$ from $\Lambda_{\varepsilon,\mu,\sigma}^\Gamma$.
Similar local problems have been studied for acoustic waves \cite{stefanov2016stable} and elastic waves \cite{stefanov2017local,bhattacharyya2018local,zhai2025determination}.
The main result of this paper is as follows.

\begin{theorem}
The three parameters $\varepsilon,\mu,\sigma$ can be uniquely determined in $\kappa^{-1}([0,1])\cap\overline{\Omega}$ by $\Lambda_{\varepsilon,\mu,\sigma}^\Gamma$, if $T$ is greater than the length of all geodesics, in the metric $c^{-2}\mathrm{d}s^2$, completely contained in $\kappa^{-1}([0,1])\cap\overline{\Omega}$.
\end{theorem}

Actually one only needs the microlocal information of $\Lambda$, that is, knowing the kernel $\Lambda(t_2,x_2,t_1,x_1)$ of $\Lambda$, up to a smooth function, on $(0,T)\times\Gamma\times(0,T)\times\Gamma$ is enough.

\section{Geometric optics}
Generally speaking, we will reduce the inverse boundary value problem to certain geometric inverse problems.
For this purpose, we would like to construct geometric optics solutions of the form
\begin{equation}\label{GOsolutions}
\mathcal{H}(t,x)\sim e^{\mathrm{i}\varrho(-t+\varphi(x))}\sum_{m=0}^\infty (\mathrm{i}\varrho)^{-m}H_m(x),\quad \mathcal{E}(t,x)=e^{\mathrm{i}\varrho(-t+\varphi(x))}\sum_{m=0}^\infty (\mathrm{i}\varrho)^{-m}E_m(x),
\end{equation}
to the equation \eqref{maineq}, with a large parameter $\varrho$. The form of the geometric optics is usually valid locally but not globally.

Let $x^1,x^2,x^3$ be a curvilinear coordinate system for $\mathbb{R}^3$ and $\mathrm{d}s^2=g_{jk}\mathrm{d}x^j\mathrm{d}x^k$ is the Euclidean metric expressed in this curvilinear coordinate. We denote $(g^{jk})=(g_{jk})^{-1}$ and $g=\det(g_{jk})$. If $\mathbf{a}=(a_1,a_2,a_3)$ is a vector field in $\mathbb{R}^3$ under the curvilinear coordinate, then
\[
(\nabla\times \mathbf{a})^i=e^{ijk}a_{k;j},
\]
where $(e^{ijk})$ is the discriminant tensor which is skew-symmetric in each pair of indices and $e^{123}=g^{-1/2}$. If in addition $\mathbf{b}=(b_1,b_2,b_3)$ is another vector field in $\mathbb{R}^3$, we have
\[
(\mathbf{a}\times\mathbf{b})^i=e^{ijk}b_ja_k.
\]


Inserting the solutions \eqref{GOsolutions} into the Maxwell's equations and collecting terms of different orders in $\varrho$, we get a system of equations
\begin{equation}\label{GOequations}
\begin{split}
\mu H_m^{~\,~i}&=e^{ijk}\nabla_j\varphi E_{mk}+e^{ijk}E_{m-1,k;j},\\
\varepsilon E_m^{~\,~i}&=-e^{ijk}\nabla_j\varphi H_{mk}-e^{ijk}H_{m-1,k;j}+\sigma E_{m-1}^{~\,~\,~\,~\,~i},
\end{split}
\end{equation}
for $m=0,1,2,\cdots$. Here we assumed that $E_{-1}=H_{-1}=0$.\\

We first consider the equation \eqref{GOequations} with $m=0$, that is,
\begin{equation}\label{eqzerothorder}
\mu H^{~i}_0=e^{ijk}\nabla_j\varphi E_{0k},\quad \varepsilon E^{~i}_0=-e^{ijk}\nabla_j\varphi H_{0k}.
\end{equation}
The above two equations can be combined into one equation for $H_0$,
\[
\mu H^{~i}_0=e^{ijk}\nabla_j\varphi \frac{1}{\varepsilon}(-e_{kpq})\nabla^p\varphi H^{~q}_0,
\]
which can be simplified to
\[
\varepsilon\mu H^{~i}_0=(-\delta_p^i\delta_q^j+\delta^j_p\delta^i_q)\nabla_j\varphi\nabla^p\varphi H_0^{~q}=|\nabla\varphi|^2 H^{~i}_0-\nabla_k\varphi H^{~k}_0\nabla^i\varphi.
\]
Here we have used the identity
\[
e^{ijk}e_{jpq}=\delta_p^i\delta_q^j-\delta^j_p\delta^i_q,
\]
which can be easily verified.
Recall that $c^2=\frac{1}{\varepsilon\mu}$. The above equation can be written as
\[
(c^{-2}-|\nabla\varphi|^2)H_0+\langle\nabla\varphi, H_0\rangle \nabla\varphi=0.
\]
Taking scalar product of the last equality with $\nabla\varphi$, we obtain
\begin{equation}\label{eqH0}
\langle\nabla\varphi, H_0\rangle=0.
\end{equation}
So we need to take
\begin{equation}\label{eikonal0}
|\nabla\varphi|^2=c^{-2}.
\end{equation}
Note also
\begin{eqnarray}
\langle E_0,\nabla\varphi\rangle&=-\frac{1}{\varepsilon}e^{ijk}\nabla_j\varphi \nabla_i\varphi H_{0k}=0,\label{eqE0}\\
\langle E_0,H_0\rangle&=-\frac{1}{\varepsilon}e^{ijk}\nabla_j\varphi H_{0i}H_{0k}=0,
\end{eqnarray}
which, together with \eqref{eqH0}, means that $\nabla\varphi,H_0,E_0$ are mutually orthogonal.

\section{boundary determination}
We will show that the parameters $\varepsilon,\mu,\sigma$, as well as all their derivatives, at $\Gamma$ can be uniquely determined by the local impedance map.
Take $x_0\in\Gamma$, and assume that $(x^1,x^2)$ are local boundary coordinates on $\Gamma$ near $x_0$.
We will work under the boundary normal coordinates specified in the following lemma. See \cite{stefanov2018inverse} for more details.
\begin{lemma}
There exists a neighborhood $N$ of $x_0$ in $\Omega$, and a diffeomorphism $\Psi: \Gamma\cap N\times [0,T)\rightarrow \Omega$ such that
\begin{enumerate}
\item $\Psi(x',0)=x'$ for all $x'\in \Gamma\cap N$;
\item $\Psi(x',x^3)=x'-x^3\nu$, where $\nu$ is the unit outer normal at the boundary point $x'$.
\end{enumerate}
In the coordinate system $(x^1,x^2,\tau=x^3)$, the Euclidean metric tensor $g_E$ takes the form
\[
g=g_{\alpha\beta}\mathrm{d}x^\alpha\mathrm{d}x^\beta+\mathrm{d}\tau\mathrm{d}\tau.
\]
\end{lemma}

From now on the Greek indices range over the values $1,2$. 
Denote
\[
J^2:=\det(g_{\alpha\beta}).
\]
Then $g=\det(g_{jk})=\det(g_{\alpha\beta})=J^2$.

We will use geometric optics solutions of the form \eqref{GOsolutions}. In $\Omega$ near $x_0$, the phase function $\varphi$ satisfies the eikonal equation \eqref{eikonal0}, which can be written as
\begin{equation}\label{varphiequation}
g^{\alpha\beta}\partial_\alpha\varphi\partial_\beta\varphi+(\partial_3\varphi)^2=c^{-2}=\varepsilon\mu,
\end{equation}
under the boundary normal coordinates.
We take
\begin{equation}\label{varpboundarycondition}
\varphi\vert_{x^3=0}=x'\cdot\xi'=x^1\xi_1+x^2\xi_2.
\end{equation}
With the extra condition $\partial_\nu\varphi\vert_{\Gamma}<0$, the equation \eqref{varphiequation} is uniquely solvable with the boundary condition \eqref{varpboundarycondition}. We take $\xi'$ such that $c^2(x)|\xi'|^2<1$ for $x$ in a neighborhood of $x_0$, so 
\[
\xi_3(x',\xi')=\partial_3\varphi(x',0)=\sqrt{c^{-2}-g^{\alpha\beta}\xi_\alpha\xi_\beta}>0.
\]
Note that $\partial_\alpha\varphi=\xi_\alpha$.

We will take 
\[
f(t,x')=e^{\mathrm{i}\varrho(-t+ x'\cdot \xi')}\chi(x',\xi'),
\]
where $\chi$ is a smooth cutoff function with small enough support near $x_0$ in $\Gamma$, and construct geometric optics solutions with
\[
\nu\times \mathcal{E}\vert_{\partial\Omega}=f.
\]
Then we have
\begin{equation}\label{equationE0boundary}
\nu\times E_0\vert_{\partial\Omega}=\chi,
\end{equation}
and
\[
\nu\times E_m\vert_{\partial\Omega}=0,\quad m=1,2,3,\cdots.
\]
Note $\nu_\alpha=0$, $\nu_3=-1$.  
We have
\[
(\nu\times H_0)^1=g^{-1/2}H_{02},\quad(\nu\times H_0)^2=-g^{-1/2}H_{01}.
\]

We write \eqref{eqzerothorder} as the following equations on $\Gamma$,
\begin{equation}\label{equationH_0boundary}
\begin{split}
\mu H^{~1}_0=e^{123}\partial_2\varphi E_{03}+e^{132}\partial_3\varphi E_{02}=g^{-1/2}\xi_2E_{03}-g^{-1/2}\xi_3 E_{02},\\
\mu H^{~2}_0=e^{213}\partial_1\varphi E_{03}+e^{231}\partial_3\varphi E_{01}=g^{-1/2}\xi_3 E_{01}-g^{-1/2}\xi_1E_{03},\\
\mu H^{~3}_0=e^{312}\partial_1\varphi E_{02}+e^{321}\partial_2\varphi E_{01}=g^{-1/2}\xi_1E_{02}-g^{-1/2}\xi_2 E_{01},
\end{split}
\end{equation}
and
\begin{equation}\label{equationE_0boundary1}
\begin{split}
\varepsilon E^{~1}_0=-e^{123}\partial_2\varphi H_{03}-e^{132}\partial_3\varphi H_{02}=-g^{-1/2}\xi_2H_{03}+g^{-1/2}\xi_3 H_{02},\\
\varepsilon E^{~2}_0=-e^{213}\partial_1\varphi H_{03}-e^{231}\partial_3\varphi H_{01}=-g^{-1/2}\xi_3 H_{01}+g^{-1/2}\xi_1H_{03},\\
\varepsilon E^{~3}_0=-e^{312}\partial_1\varphi H_{02}-e^{321}\partial_2\varphi H_{01}=-g^{-1/2}\xi_1H_{02}+g^{-1/2}\xi_2 H_{01}.
\end{split}
\end{equation}
Using the equalities
\[
g^{11}=J^{-2}g_{22},\quad g^{22}=J^{-2}g_{11},\quad g^{12}=-J^{-2}g_{12},
\]
we calculate
\[
\begin{split}
(\nu\times E_0)_1&=g_{11}e^{132}\nu_3E_{02}+g_{12}e^{231}\nu_3E_{01}\\
&=g^{-1/2}(g_{11}E_{02}-g_{21}E_{01})\\
&=g^{-1/2}(g_{11}g_{21}E_0^{~1}+g_{11}g_{22}E_0^{~2}-g_{21}g_{11}E_0^{~1}-g_{21}g_{12}E_0^{~2})\\
&=\varepsilon^{-1}g^{-1}J^2\Big(g_{11}g^{12}(\xi_2H_{03}-\xi_3H_{02})+g_{11}g^{11}(-\xi_3H_{01}+\xi_1H_{03})\\
&\quad\quad\quad+g_{12}g^{22}(\xi_2H_{03}-\xi_3H_{02})+g_{21}g^{12}(-\xi_3H_{01}+\xi_1H_{03})\Big)\\
&=-\varepsilon^{-1}(\xi_3H_{01}-\xi_1H_{03}),
\end{split}
\]
and similarly,
\[
\begin{split}
(\nu\times E_0)_2
=-\varepsilon^{-1}(\xi_3H_{02}-\xi_2H_{03}).
\end{split}
\]
Note that we have the identity
\begin{equation}\label{formulaH_03}
H_{03}=-\frac{1}{\xi_3}g^{\alpha\beta}H_{0\alpha}\xi_\beta,
\end{equation}
resulting from \eqref{equationH_0boundary}.\\

Take $\chi$ such that $\chi=\xi'=(\xi_1,\xi_2,0)$ in a neighborhood of $x_0$. Then one can solve the equation \eqref{equationE0boundary}, which can be written as
\[
-\varepsilon^{-1}(\xi_3H_{0\alpha}+\frac{1}{\xi_3}\xi^\beta H_{0\beta}\xi_\alpha)=\xi_\alpha,
\]
 for $H_0$ to obtain
\begin{equation}\label{formulaH_0}
H_{0\alpha}=-\varepsilon \left(\xi_3+\frac{|\xi'|^2}{\xi_3}\right)^{-1}\xi_\alpha.
\end{equation}

So from the impedance map we can recover
\[
\begin{split}
(\xi'\times (\nu\times H_0))^3=&e^{312}\xi_1(\nu\times H_0)_2+e^{321}\xi_2(\nu\times H_0)_1\\
=&g^{-1/2}\xi_1(-g_{21}g^{-1/2}H_{02}+g_{22}g^{-1/2}H_{01})-g^{-1/2}\xi_2(-g_{11}g^{-1/2}H_{02}+g_{12}g^{-1/2}H_{01})\\
=&g^{-1}(-g_{21}\xi_1H_{02}+g_{22}\xi_1H_{01}+g_{11}\xi_2H_{02}-g_{12}\xi_2H_{01})\\
=&(g^{12}\xi_1H_{02}+g^{11}\xi_1H_{01}+g^{22}\xi_2H_{02}+g^{12}\xi_2H_{01})\\
=&-\varepsilon\left(\xi_3+\frac{|\xi'|^2}{\xi_3}\right)^{-1}|\xi'|^2
\end{split}
\]
in a neighborhood of $x_0$.
Recall that
\[
\xi_3(x,\xi')=\sqrt{\varepsilon\mu-|\xi'|^2}.
\]
We can recover
\[
(\xi'\times (\nu\times H_0))^3=-\frac{1}{\mu(x)}\sqrt{\varepsilon(x)\mu(x)-|\xi'|^2}|\xi'|^2.
\]
Varying the values of $|\xi'|$, one can recover $\mu$ and $\varepsilon$ on $\Gamma$ since $x_0$ can be any point on $\Gamma$. Using equations \eqref{equationH_0boundary}, \eqref{equationE_0boundary1}, we can recover $E_0,H_0$ on $\Gamma$.\\

Next, we recover the normal derivatives $\partial_3\mu$, $\partial_3\varepsilon$ and $\sigma$ on $\Gamma$.
Notice that the equation
\[
\nu\times E_1=0,
\]
implies $E_{1\alpha}=0$ for $\alpha=1,2$.

First we take derivatives of the eikonal equation \eqref{varphiequation} in $x^3$, we have
\[
2\partial_3\varphi\partial^2_{33}\varphi=\varepsilon\partial_3\mu+\mu\partial_3\varepsilon+R.
\]
Here and in the section $R$ denotes terms that depend on $\mu\vert_{\Gamma},\varepsilon\vert_{\Gamma}$.
So
\begin{equation}\label{secondderivativephi}
\partial^2_{33}\varphi=\frac{1}{2\xi_3}(\varepsilon\partial_3\mu+\mu\partial_3\varepsilon)+R.
\end{equation}

Taking derivatives $\frac{\partial}{\partial x^3}$ of the equation \eqref{eqzerothorder}, we obtain equations
\[
\partial_3\mu H_0^{~1}+\mu\partial_3H_0^{~1}=g^{-1/2}\xi_2\partial_3E_{03}-g^{-1/2}\partial^2_{3}\varphi E_{02}-g^{-1/2}\xi_3\partial_3E_{02}+R,
\]
\[
\partial_3\mu H_0^{~2}+\mu\partial_3H_0^{~2}=-g^{-1/2}\xi_1\partial_3E_{03}+g^{-1/2}\partial^2_{3}\varphi E_{01}+g^{-1/2}\xi_3\partial_3E_{01}+R,
\]
on $\Gamma$.
Summing the above two equations, we have
\begin{equation}\label{eq123}
\begin{split}
&-g^{-1/2}\xi_3(\xi_1\partial_3E_{02}-\xi_2\partial_3E_{01})\\
=&\partial_3\mu(\xi_1H_0^{~1}+\xi_2H_0^{~2})+\mu(\xi_1\partial_3H_0^{~1}+\xi_2\partial_3H_0^{~2})+g^{-1/2}\xi_1\partial^2_{3}\varphi E_{02}-g^{-1/2}\xi_2\partial^2_{3}\varphi E_{01}+R.
\end{split}
\end{equation}

Note that by \eqref{formulaH_0} and \eqref{formulaH_03},
\[
\xi_1H_0^{~1}+\xi_2H_0^{~2}=-\frac{1}{\mu}\xi_3|\xi'|^2,
\]
and
\[
H_{03}=\frac{1}{\mu}|\xi'|^2.
\]
Using \eqref{equationE_0boundary1}, we have
\[
\varepsilon(\xi^2E_0^{~1}-\xi^1E_0^{~2})=-g^{-1/2}|\xi'|^2H_{03}+g^{-1/2}\xi_3(\xi^2H_{02}+\xi^1H_{01})=-\varepsilon g^{-1/2}|\xi'|^2.
\]
So
\[
\xi_2E_{01}-\xi_1E_{02}=\xi^2E_0^{~1}-\xi^1E_0^{~2}=-g^{-1/2}|\xi'|^2.
\]

Taking $m=1$ in the equation \eqref{GOequations}, we have
\[
\mu H_1^{~1}=g^{-1/2}\xi_2 E_{13}-g^{-1/2}\partial_3E_{02}+R,
\]
\[
\mu H_1^{~2}=-g^{-1/2}\xi_1E_{13}+g^{-1/2}\partial_3E_{01}+R.
\]
Then
\begin{equation}\label{eqH1comb}
\mu(\xi_1H^{~1}_1+\xi_2H_1^{~2})=-g^{-1/2}\xi_1\partial_3E_{02}+g^{-1/2}\xi_2\partial_3E_{01}+R.
\end{equation}
We also have
\[
0=\varepsilon E_1^{~1}=-g^{-1/2}\xi_2H_{13}+g^{-1/2}\xi_3H_{12}+g^{-1/2}\partial_3H_{02}+\sigma E_0^{~1}+R,
\]
\[
0=\varepsilon E_1^{~2}=g^{-1/2}\xi_1H_{13}-g^{-1/2}\xi_3H_{11}-g^{-1/2}\partial_3H_{01}+\sigma E_0^{~2}+R,
\]
and
\[
\mu H_1^{~3}=R.
\]
So
\[
\partial_3H_{02}=-\xi_3H_{12}-\sigma g^{1/2}E_0^{~1}+R,\quad \partial_3H_{01}=-\xi_3H_{11}+\sigma g^{1/2}E_0^{~2}+R.
\]

Using also the equation \eqref{eqH1comb}, we have
\[
\begin{split}
\xi_1\partial_3H_{0}^{~1}+\xi_2\partial_3H_{0}^{~2}=&\xi^1\partial_3H_{01}+\xi^2\partial_3H_{02}+R\\
=&-\xi_3(\xi^2H_{12}+\xi^1H_{11})+\sigma g^{1/2}\xi^1E_0^{~2}-\sigma g^{1/2}\xi^2E_0^{~1}+R\\
=&\frac{1}{\mu}g^{-1/2}\xi_3(\xi_1\partial_3E_{02}-\xi_2\partial_3E_{01})+\sigma g^{1/2}\xi^1E_0^{~2}-\sigma g^{1/2}\xi^2E_0^{~1}+R.
\end{split}
\]
Inserting into \eqref{eq123}, and using \eqref{secondderivativephi}, we have
\[
\begin{split}
&2\mu(\xi_1\partial_3H_0^{~1}+\xi_2\partial_3H_0^{~2})\\
=&-\partial_3\mu(\xi_1H_0^{~1}+\xi_2H_0^{~2})-\xi_1\partial^2_{3}\varphi E_{02}+\xi_2\partial^2_{3}\varphi E_{01}+\sigma\mu g^{1/2}\xi^1E_0^{~2}-\sigma\mu g^{1/2}\xi^2E_0^{~1}+R\\
=&\frac{|\xi'|^2}{\xi_3}\left(\frac{\xi_3^2}{\mu}\partial_3\mu-\frac{1}{2}(\varepsilon\partial_3\mu+\mu\partial_3\varepsilon)\right)+\sigma\mu |\xi'|^2+R\\
=&\frac{|\xi'|^2}{\sqrt{\varepsilon\mu-|\xi'|^2}}\left((\varepsilon\mu-|\xi'|^2)\frac{\partial_3\mu}{\mu}-\frac{1}{2}((\varepsilon\mu)\frac{\partial_3\mu}{\mu}+\mu\partial_3\varepsilon)\right)+\sigma\mu |\xi'|^2+R.
\end{split}
\]

From the impedance map, we can recover
\[
\begin{split}
(\xi'\times (\nu\times H_1))^3=\xi^1H_{11}+\xi^2H_{12}=&-\frac{1}{\xi_3}(\xi^1\partial_3H_{01}+\xi^2\partial_3H_{02})+R\\
=&-\frac{1}{\xi_3}(\xi_1\partial_3H_{0}^1+\xi_2\partial_3H_{0}^2)+R
\end{split}
\]
in a neighborhood of $x_0$, and thus can recover 
\[
(\varepsilon\mu-|\xi'|^2)\frac{\partial_3\mu}{\mu}-\frac{1}{2}((\varepsilon\mu)\frac{\partial_3\mu}{\mu}+\mu\partial_3\varepsilon)+\sigma \mu\sqrt{\varepsilon\mu-|\xi'|^2}
\]
on $\Gamma$.
By varying the values of $|\xi'|$, we can recover $\partial_3\mu\vert_{\Gamma}$, $\partial_3\varepsilon\vert_{\Gamma}$ and $\sigma\vert_\Gamma$ separately. Note that now we can also get the values of $\partial_3E_0,\partial_3H_0$ and $E_{1},H_1$ on $\Gamma$.

Repeating this argument we can recover $\partial^m_3\mu\vert_{\Gamma}$  $\partial^m_3\mu\vert_{\Gamma}$ and $\partial^{m-1}_3\sigma\vert_{\Gamma}$ for all $m\geq 2$. We summarize the results of this section in the following theorem.
\begin{theorem}\label{boundaryjets}
The local impedance map $\Lambda_{\lambda,\mu,\sigma}^\Gamma$ uniquely  determines $\partial^m_\nu\varepsilon,\partial^m_\nu\mu,\partial^m_\nu\sigma$, $m=0,1,2,\cdots$, on $\Gamma$.
\end{theorem}

\section{Interior determination}
We will consider the Cauchy problem for the Maxwell's equations
\begin{equation}
\begin{cases}
\mu\partial_t\mathcal{H}=-\nabla\times \mathcal{E},\quad \varepsilon\partial_t\mathcal{E}=\nabla\times \mathcal{H}+\sigma\mathcal{E}\quad \text{in }\mathbb{R}^3\times (0,T),\\
\mathcal{H}(0,\cdot)=h,\quad \mathcal{E}(0,\cdot)=e\quad \quad\quad\quad\text{in }\mathbb{R}^3,
\end{cases}
\end{equation}
and construct solutions of the form
\begin{equation}\label{FIOinterior}
\mathcal{H}(t,x)=(2\pi)^{-3}\int_{\mathbb{R}^3} e^{\mathrm{i}\phi(t,x,\xi)}H(t,x,\xi)\hat{h}(\xi)\mathrm{d}\xi,\quad \mathcal{E}(t,x)=(2\pi)^{-3}\int_{\mathbb{R}^3} e^{\mathrm{i}\phi(t,x,\xi)}E(t,x,\xi)\hat{e}(\xi)\mathrm{d}\xi,
\end{equation}
where
\[
H(t,x,\xi)\sim\sum_{m=0}^\infty \mathrm{i}^{-m}H_m(t,x,\xi),\quad E(t,x,\xi)\sim\sum_{m=0}^\infty \mathrm{i}^{-m}E_m(t,x,\xi),
\]
and $H_m,E_m$ are homogeneous in $\xi$ of order $-m$ for large $|\xi|$.

The phase function satisfies the eikonal equation
\begin{equation}
(\partial_t\phi)^2=c^2|\nabla_x\phi|^2,
\end{equation}
with initial value
\[
\phi\vert_{t=0}=x\cdot\xi.
\]

The integrals in \eqref{FIOinterior} are Fourier integral operators (FIOs) \cite{hormander1971fourier,duistermaat1994fourier}. The singularities of $(h,e)$ propagate along the null bicharacteristics of the form $(t(s),x(s),\tau(s),\xi(s))$ where
\[
\frac{\mathrm{d}t}{\mathrm{d}s}=1,\quad \frac{\mathrm{d}x}{\mathrm{d}s}=\pm c\frac{\xi}{|\xi|},\quad \frac{\mathrm{d}\tau}{\mathrm{d}s}=0,\quad \frac{\mathrm{d}\xi}{\mathrm{d}s}=\mp(\nabla_xc)|\xi|.
\]
Note that $(x(s),\xi(s))$ is then the geodesics with respect to the metric $g=c^{-2}g_E$ lifted to the phase space. We can write
\[
\mathrm{WF}(\mathcal{H},\mathcal{E})=C_+\circ\mathrm{WF}(h,e)\cup C_-\circ\mathrm{WF}(h,e),
\]
where
\[
\begin{split}
C_+(x,\xi)&=(s,\gamma_{x,\xi/|\xi|_g}(s),-|\xi|_g,|\xi|_gg\dot{\gamma}_{x,\xi/|\xi|_g}(s)),\\
C_-(x,\xi)&=(s,\gamma_{x,-\xi/|\xi|_g}(s),|\xi|_g,-|\xi|_gg\dot{\gamma}_{x,-\xi/|\xi|_g}(s)).
\end{split}
\]
Here $\gamma_{x,\eta}$ is the geodesic issued from $x$ in the direction $g^{-1}\eta$.
If we consider $t$ as a parameter, we have
\[
\mathrm{WF}(\mathcal{H}(t,\cdot),\mathcal{E}(t,\cdot))=C_+(t)\circ\mathrm{WF}(h,e)\cup C_-(t)\circ\mathrm{WF}(h,e),
\]
where
\[
\begin{split}
C_+(t)(x,\xi)&=(\gamma_{x,\xi/|\xi|_g}(s),|\xi|_gg\dot{\gamma}_{x,\xi/|\xi|_g}(s)),\\
C_-(t)(x,\xi)&=(\gamma_{x,-\xi/|\xi|_g}(s),-|\xi|_gg\dot{\gamma}_{x,-\xi/|\xi|_g}(s)).
\end{split}
\]
The above construction can be done in some neighborhood of a fixed point $(0,x_0)$ in general. To extend globally, we can first localize it first for $(h,e)$ with $WF((h,e))$ in a neighborhood of some fixed $(x_0,\xi_0)\in T\mathbb{R}^3\setminus 0$. Then $(\mathcal{H},\mathcal{E})$ will be defined near the geodesic issued from $(x_0,\xi_0)$ but in some neighborhood of $(0,x_0)$. Then we can fix some $t=t_1>0$ at which the solution is still defined, replace $\{t=0\}$ with $\{t=t_1\}$, take the Cauchy data there, and then construct a new FIO solution of the same form. Continuing with this process, we can make the construction global.

Since the jets of $\varepsilon,\mu,\sigma$ on $\Gamma$ can be determined (cf. Theorem \ref{boundaryjets}), we can extend $\varepsilon,\mu,\sigma$ smoothly to some small neighborhood $U$ of $\Gamma$. Choose $\zeta_1=(t_1,x_1,\tau_1=-c|\xi_1|,\xi_1)\in T^*((0,\epsilon)\times\overline{\Omega})\setminus 0$, with $x_1\in\Gamma$ and $\xi_1$ pointing into $\Omega$. Assume that the null bicharacteristic $\vartheta$ passing though $\zeta_1$ is transversal at that point, and hits $T^*(R\times\Gamma)$ again transversely at some $\zeta_2=(t_2,x_2,\tau_2=-c|\xi_2|,\xi_2)$ with $x_2\in\Gamma$. Extend $\vartheta$ outside the domain $\Omega$ until it hits $\{t=0\}$ at some point $\zeta_0=(t_0,x_0,\tau_0=-c|\xi_0|,\xi_0)$. For $\epsilon>0$ small enough, $x_0\in U\setminus\overline{\Omega}$. 

We can construct microlocal solutions $\mathcal{H}^{\mathrm{inc}},\mathcal{E}^{\mathrm{inc}}$ of the form \eqref{FIOinterior} propagating along $\vartheta$. We can cut $\mathcal{H}^{\mathrm{inc}},\mathcal{E}^{\mathrm{inc}}$ smoothly so that its support is concentrated near $\vartheta$. Then $(\mathcal{H}^{\mathrm{inc}}(0,\cdot),\mathcal{E}^{\mathrm{inc}}(0,\cdot))=(h,e)=0$ in $\Omega$. The trace of $\nu\times \mathcal{E}^{\mathrm{inc}}$ on $\Gamma$ can be written as $f_1+f_2$, where $f_j$ is supported in a small neighborhood of $(t_j,x_j)$. Let $\mathcal{H},\mathcal{E}$ be the solutions to \eqref{maineq}, \eqref{initialcondition}, \eqref{boundarycondition} with $f=f_1$. Then $(\mathcal{H},\mathcal{E})$ and $(\mathcal{H}^{\mathrm{inc}},\mathcal{E}^{\mathrm{inc}})$ differ by a smooth function before they hit $\partial\Omega$ again. To satisfy the zero boundary condition $\nu\times\mathcal{E}=0$ near $(t_2,x_2)$, we can write $(\mathcal{H},\mathcal{E})$ as the sum of the incident wave $(\mathcal{H}^{\mathrm{inc}},\mathcal{E}^{\mathrm{inc}})$ and the reflected wave $(\mathcal{H}^{\mathrm{ref}},\mathcal{E}^{\mathrm{ref}})$: 
\[
(\mathcal{H},\mathcal{E})=(\mathcal{H}^{\mathrm{inc}},\mathcal{E}^{\mathrm{inc}})+(\mathcal{H}^{\mathrm{ref}},\mathcal{E}^{\mathrm{ref}}),
\]
where
\[
\begin{split}
\mathcal{H}^{\mathrm{inc}}(t,x)&=(2\pi)^{-3}\int_{\mathbb{R}^3} e^{\mathrm{i}\phi(t,x,\xi)}(\sum_{m=0}^\infty \mathrm{i}^{-m}H_m^{\mathrm{inc}}(t,x,\xi))\hat{h}(\xi)\mathrm{d}\xi,\\
\mathcal{E}^{\mathrm{inc}}(t,x)&=(2\pi)^{-3}\int_{\mathbb{R}^3} e^{\mathrm{i}\phi(t,x,\xi)}(\sum_{m=0}^\infty \mathrm{i}^{-m}E_m^{\mathrm{inc}}(t,x,\xi))\hat{e}(\xi)\mathrm{d}\xi,\\
\mathcal{H}^{\mathrm{ref}}(t,x)&=(2\pi)^{-3}\int_{\mathbb{R}^3} e^{\mathrm{i}\phi^{\mathrm{ref}}(t,x,\xi)}(\sum_{m=0}^\infty \mathrm{i}^{-m}H_m^{\mathrm{ref}}(t,x,\xi))\hat{h}(\xi)\mathrm{d}\xi,\\
\mathcal{E}^{\mathrm{ref}}(t,x)&=(2\pi)^{-3}\int_{\mathbb{R}^3} e^{\mathrm{i}\phi^{\mathrm{ref}}(t,x,\xi)}(\sum_{m=0}^\infty \mathrm{i}^{-m}E_m^{\mathrm{ref}}(t,x,\xi))\hat{e}(\xi)\mathrm{d}\xi.
\end{split}
\]
Here the phase function $\phi^{\mathrm{ref}}$ solves the same eikonal equation as $\phi$ but satisfies the boundary condition $\phi^{\mathrm{ref}}=\phi$ near $(t_2,x_2)$. It differs from $\phi$ by the sign of its normal derivative $\partial_\nu\phi=-\partial_\nu\phi^{\mathrm{ref}}$ near $(t_2,x_2)$. So, under the boundary normal coordinates,
\[
E^{\mathrm{ref}}_{m\alpha}=-E^{\mathrm{inc}}_{m\alpha}
\]
for $\alpha=1,2$, $m=0,1,2,\cdots$. Note that
\[
\nu\times\mathcal{H}=(2\pi)^{-3}\int_{\mathbb{R}^3} e^{\mathrm{i}\phi(t,x,\xi)}(\sum_{m=0}^\infty \nu\times(H_m^{\mathrm{inc}}+H_m^{\mathrm{ref}})(t,x,\xi))\hat{h}(\xi)\mathrm{d}\xi.
\]
By the third equation in \eqref{equationH_0boundary}, we have
\[
\mu (H_0^{3,\mathrm{inc}}+H_0^{3,\mathrm{ref}})=-g^{-1/2}\xi_1(E_{02}^\mathrm{inc}+E_{02}^\mathrm{ref})+g^{-1/2}\xi_2(E_{01}^\mathrm{inc}+E_{01}^\mathrm{ref})=0.
\]
By the first two equations in \eqref{equationE_0boundary1}, using the fact $\partial_3\phi^\mathrm{ref}=-\partial_3\phi$, we have
\[
\begin{split}
\xi_2( H_{03}^\mathrm{inc}+ H_{03}^\mathrm{ref})-\xi_3(H_{02}^\mathrm{inc}-H_{02}^\mathrm{ref})=0,\\
\xi_1( H_{03}^\mathrm{inc}+ H_{03}^\mathrm{ref})-\xi_3(H_{01}^\mathrm{inc}-H_{01}^\mathrm{ref})=0.
\end{split}
\]
Therefore we obtain
\[
H_{0\alpha}^\mathrm{inc}=H_{0\alpha}^\mathrm{ref},\quad \alpha=1,2.
\]
From $\nu\times\mathcal{H}\vert_{\Gamma}$ we can recover $H_{0\alpha}^\mathrm{inc}+H_{0\alpha}^\mathrm{ref}$ and thus can recover $H_{0\alpha}^\mathrm{inc}$ on $\Gamma$. We can then get $E_{03}^\mathrm{inc}$ using the last equation in \eqref{equationE_0boundary1}. The using the quations in \eqref{equationH_0boundary}, we obtain $E_{0\alpha}^\mathrm{inc}$ and subsequently $H_{03}^\mathrm{inc}$. Continuing with this process, we can recover $H_{m}^\mathrm{inc}\vert_{(0,T)\times\Gamma}$, $E_{m}^\mathrm{inc}\vert_{(0,T)\times\Gamma}$ for all $m=1,2,\cdots$. Therefore, we can recover the full microlocal solution $(\mathcal{H}^{\mathrm{inc}},\mathcal{E}^{\mathrm{inc}})\vert_{(0,T)\times\Gamma}$ from the impedance map. So from now on, we ignore the boundary reflection and use the notation $(\mathcal{H},\mathcal{E})$ in place of $(\mathcal{H}^{\mathrm{inc}},\mathcal{E}^{\mathrm{inc}})$.

\subsection{Recovery of the wavespeed}
Summarize above discussions, the impedance map determines the wavefront set of $(\mathcal{H}^{\mathrm{inc}},\mathcal{E}^{\mathrm{inc}})\vert_{(0,T)\times\Gamma}$, $WF((\mathcal{H}^{\mathrm{inc}},\mathcal{E}^{\mathrm{inc}})\vert_{(0,T)\times\Gamma})$. Since the wave speed $c$ is already known on $\Gamma$ by the boundary determination results, we indeed have $WF((\mathcal{H}^{\mathrm{inc}},\mathcal{E}^{\mathrm{inc}}))\vert_{(0,T)\times\Gamma}$. This in turn determines the \textit{lens relations} on $\Gamma$, which we shall briefly review now. We refer to \cite{stefanov2018inverse,stefanov2017local} for more detailed discussions. Consider the Riemannian manifold $(\Omega,g)$ with smooth boundary $\partial\Omega$. Denote
\[
\partial_\pm S\Omega:=\{(x,v)\in T\Omega, x\in\partial\Omega,|v|_g=1, \pm\langle v,\nu\rangle_g\leq 0\},
\]
where $\nu$ is the unit outer normal vector to $\partial\Omega$. Assume that $(\Omega,g)$ is non-trapping. For any $(x,v)\in\partial_+S\Omega$, let $\gamma_{x,v}$ be the unique geodesic determined by $(x,v)$ such that $\gamma_{x,v}(0)=0$ and $\dot{\gamma}_{x,v}(0)=v$. Define the \textit{exist time} $\tau(x,v)\in [0,\infty)$ so that
\[
\tau(x,v)=\sup\{t,\gamma_{x,v}(t)\in\Omega\}.
\]
We define the \textit{scattering relation}
\[
\alpha:\partial_+S\Omega\rightarrow \partial_-S\Omega,
\]
in the way such that $\alpha(x,v)=(y,w)$, where $y=\gamma_{x,v}(\tau(x,v))$ and $w=\dot{\gamma}_{x,v}(\tau(x,v))$. The maps $(\tau,\alpha)$ together are called the \textit{lens relation}.

By the above discussion, the local impedance map determines the lens relation $(\tau(x,v),\alpha(x,v))$ for $x\in\Gamma$ and $\gamma_{x,v}$ hits $\Gamma$ again transversely before time $T$.
Therefore, using the results in \cite{stefanov2016boundary}, we can recover the wavespeed $c$, or equivalently, $\varepsilon\mu$, in $\kappa^{-1}([0,1])\cap\overline{\Omega}$.

\subsection{Recovery of $\frac{\sigma}{\varepsilon}$}
To recover $\frac{\sigma}{\varepsilon}$ with $c$ already known, we need to have an explicit construction of the amplitudes. 
We will use solutions of the form
\[
\mathcal{H}(t,x)=(2\pi)^{-1}\int_\mathbb{R} e^{\mathrm{i}\varrho(-t+\varphi(x))}H(x,\varrho)\mathrm{d}\varrho,\quad \mathcal{E}(t,x)=(2\pi)^{-1}\int_\mathbb{R} e^{\mathrm{i}\varrho(-t+\varphi(x))}E(x,\varrho)\mathrm{d}\varrho,
\]
where
\[
H\sim\sum_{m=0}^\infty (\mathrm{i}\varrho)^{-m}H_m(x),\quad E\sim\sum_{m=0}^\infty (\mathrm{i}\varrho)^{-m}E_m(x).
\]
Here the phase function $\varphi$ satisfies the eikonal equation \eqref{eikonal0}. This amounts to a reparametrization of the Langragian submanifold associated with \eqref{FIOinterior}.

We will carry out the calculations in the so-called ray coordinates for $\mathbb{R}^3$, i.e., curvilinear coordinates $(x^1,x^2,x^3)$ such that $x^3(x)=\tau=\varphi(x)$ and the coordinate surface $x^3=x^3_0$ are orthogonal to the coordinate lines $x^1=x_0^1$, $x^2=x^2_0$ that are the geodesics of metric $g=c^{-2}\mathrm{d}s^2$. Under these ray coordinates the Euclidean metric has the form
\[
\mathrm{d}s^2=g_{\alpha\beta}\mathrm{d}x^\alpha\mathrm{d}x^\beta+c^2\mathrm{d}\tau^2.
\]
The Christoffel symbols of the Euclidean metric are then
\[
\begin{split}
&\Gamma^3_{\alpha\beta}=-\frac{1}{2}c^{-2}\frac{\partial g_{\alpha\beta}}{\partial\tau},\quad \Gamma^\beta_{\alpha 3}=\frac{1}{2}g^{\beta\gamma}\frac{\partial g_{\alpha\gamma}}{\partial \tau},\quad \Gamma^3_{\alpha 3}=c^{-1}\frac{\partial c}{\partial x^\alpha},\\
&\Gamma^\alpha_{33}=-cg^{\alpha\beta}\frac{\partial c}{\partial x^\beta},\quad \Gamma^3_{33}=c^{-1}\frac{\partial c}{\partial \tau}.
\end{split}
\]
We use the notation $J^2:=\det(g_{\alpha\beta})$. Then $g=\det(g_{jk})=c^2J^2$. We follow the lines of calculation in \cite{Shara} closely.

Under the above ray coordinates,
\[
\varphi_{;3}=1,\varphi_{;}^{~~3}=c^{-2},\varphi_{;\alpha}=\varphi_{;}^{~~\alpha}=0.
\]
The equation \eqref{eqzerothorder} can be written as
\[
\mu H_0^{~\alpha}=e^{\alpha3\beta}E_{0\beta},\quad \varepsilon E^{~\alpha}_0=-e^{\alpha3\beta }H_{0\beta}.
\]
Then by \eqref{eqH0} and \eqref{eqE0} we have
\[
E_{03}=E^{~3}_0=H_{03}=H^{~3}_0=0.
\]
and therefore
\begin{equation}\label{E0H0relation}
E_0^{~1}=\frac{1}{\varepsilon}g^{-1/2}H_{02},\quad E_0^{~2}=-\frac{1}{\varepsilon}g^{-1/2}H_{01}.
\end{equation}

Taking $m=1$ and putting $i=\alpha$ in \eqref{GOequations}, we have
\[
\mu H_1^{~\alpha}=e^{\alpha 3\beta }E_{1\beta}+e^{\alpha jk}E_{0k;j},
\]
\[
\varepsilon E_1^{~\alpha}=-e^{\alpha 3\beta }H_{1\beta}-e^{\alpha jk}H_{0k;j}+\sigma E_0^{~\alpha}.
\]
Taking $\alpha=1$ and $\alpha=2$, we write
\begin{align}
\mu H_1^{~1}&=-g^{-1/2}E_{12}-g^{-1/2}(E_{02;3}-E_{03;2}),\label{eqH1}\\
\mu H_1^{~2}&=g^{-1/2}E_{11}+g^{-1/2}(E_{01;3}-E_{03;1}),\label{eqH2}\\
\varepsilon E_1^{~1}&=g^{-1/2}H_{12}+g^{-1/2}(H_{02;3}-H_{03;2})+\sigma E_0^{~1},\label{eqE1}\\
\varepsilon E_1^{~2}&=-g^{-1/2}H_{11}-g^{-1/2}(H_{01;3}-H_{03;1})+\sigma E_0^{~2}.\label{eqE2}
\end{align}
Upon using $H_{1\alpha}=g_{\alpha\beta}H_{1}^{~\beta}$ in \eqref{eqH1} and \eqref{eqH2}, we have
\[
\begin{split}
\mu H_{11}=&\mu H_1^{~1}g_{11}+\mu H_1^{~2}g_{12}\\
=& -g^{-1/2}E_{12}g_{11}-g^{-1/2}(E_{02;3}-E_{03;2})g_{11} +g^{-1/2}E_{11}g_{12}+g^{-1/2}(E_{01;3}-E_{03;1})g_{12}\\
=& -g^{-1/2}E_{12}c^{-2}gg^{22}-g^{-1/2}(E_{02;3}-E_{03;2})c^{-2}gg^{22}-g^{-1/2}E_{11}c^{-2}gg^{12}-g^{-1/2}(E_{01;3}-E_{03;1})gg^{12}\\
=&-g^{1/2}c^{-2} (E_1^{~2}+E_{0~;3}^{~2}-E_{03;}^{~\,~\,~2}),
\end{split}
\]
and similarly
\[
\mu H_{12}=g^{1/2}c^{-2} (E_1^{~1}+E_{0~;3}^{~1}-E_{03;}^{~\,~\,~1}).
\]
Inserting into \eqref{eqE1} and \eqref{eqE2}, we obtain
\begin{equation}\label{eqamplitude0}
\begin{split}
c^{-2}(E_{0~;3}^{~1}-E_{03;}^{~\,~\,1})+\mu g^{-1/2}(H_{02;3}-H_{03;2})+\mu\sigma E_0^{~1}=0,\\
c^{-2}(E_{0~;3}^{~2}-E_{03;}^{~\,~\,2})-\mu g^{-1/2}(H_{01;3}-H_{03;1})+\mu\sigma E_0^{~2}=0.
\end{split}
\end{equation}
Note
\[
\begin{split}
E_{0~;3}^{~1}=&\partial_3E_0^{~1}+\Gamma^1_{3\alpha}E_0^{~\alpha}\\
=&\partial_3\left(\frac{1}{\varepsilon}g^{-1/2}H_{02}\right)+\Gamma^1_{31}\frac{1}{\varepsilon}g^{-1/2}H_{02}-\Gamma^1_{32}\frac{1}{\varepsilon}g^{-1/2}H_{01}\\
=&\partial_3g^{-1/2}\varepsilon^{-1}H_{02}+g^{-1/2}\frac{\partial (\varepsilon^{-1}H_{02})}{\partial\tau}+\frac{1}{2}g^{1\gamma}\frac{\partial g_{1\gamma}}{\partial\tau} g^{-1/2}\varepsilon^{-1}H_{02}-\frac{1}{2}g^{1\gamma}\frac{\partial g_{2\gamma}}{\partial\tau}g^{-1/2}\varepsilon^{-1}H_{01},
\end{split}
\]
\[
\begin{split}
E_{0~;3}^{~2}=&\partial_3E_0^{~2}+\Gamma^2_{3\alpha}E_0^{~\alpha}\\
=&\partial_3\left(-\frac{1}{\varepsilon}g^{-1/2}H_{01}\right)+\Gamma^2_{31}\frac{1}{\varepsilon}g^{-1/2}H_{02}-\Gamma^2_{32}\frac{1}{\varepsilon}g^{-1/2}H_{01}\\
=&-\partial_3g^{-1/2}\varepsilon^{-1}H_{01}-g^{-1/2}\frac{\partial (\varepsilon^{-1}H_{01})}{\partial\tau}+\frac{1}{2}g^{2\gamma}\frac{\partial g_{1\gamma}}{\partial\tau} g^{-1/2}\varepsilon^{-1}H_{02}-\frac{1}{2}g^{2\gamma}\frac{\partial g_{2\gamma}}{\partial\tau}g^{-1/2}\varepsilon^{-1}H_{01}.
\end{split}
\]
Next we calculate
\[
\begin{split}
E_{03;}^{~\,~\,~\alpha}=&g^{\alpha\beta}E_{03;\beta}=-g^{\alpha\beta}\Gamma_{3\beta}^\gamma E_{0}^{~\delta}g_{\delta\gamma}=-g^{\alpha\beta}\Gamma^\gamma_{3\beta}(\frac{1}{\varepsilon}g^{-1/2}g_{1\gamma}H_{02}-\frac{1}{\varepsilon}g^{-1/2}g_{2\gamma}H_{01})\\
=&-\frac{1}{2}\frac{\partial g_{1\beta}}{\partial\tau}g^{\alpha\beta}g^{-1/2}\varepsilon^{-1}H_{02}+\frac{1}{2}\frac{\partial g_{2\beta}}{\partial\tau}g^{\alpha\beta}g^{-1/2}\varepsilon^{-1}H_{01}.
\end{split}
\]

Then we can write the equations \eqref{eqamplitude0} as
\begin{equation}
\begin{split}
2\varepsilon^{-1}g^{-1/2}\frac{\partial H_{0\alpha}}{\partial\tau}+&\frac{\partial\varepsilon^{-1}}{\partial\tau}g^{-1/2}H_{0\alpha}+g^{\beta\gamma}\frac{\partial g_{\beta\gamma}}{\partial\tau}g^{-1/2}\varepsilon^{-1}H_{0\alpha}\\
&-\frac{1}{2}g^{-1/2}\frac{1}{g}\frac{\partial g}{\partial\tau}\varepsilon^{-1}H_{0\alpha}-\frac{\partial g_{\alpha\beta}}{\partial\tau}g^{-1/2}\varepsilon^{-1}H_0^{~\beta}+\frac{\sigma}{\varepsilon^2}g^{-1/2}H_{0\alpha}=0.
\end{split}
\end{equation}
Using the identity
\[
g^{\alpha\beta}\frac{\partial g_{\alpha\beta}}{\partial\tau}=\frac{1}{g}\frac{\partial g}{\partial\tau}-\frac{2}{c}\frac{\partial c}{\partial \tau}=\frac{2}{J}\frac{\partial J}{\partial\tau},
\]
the above equation simplifies to 
\[
2\frac{\partial H_{0\alpha}}{\partial\tau}-\frac{1}{\varepsilon} \frac{\partial\varepsilon}{\partial\tau}H_{0\alpha}+\frac{1}{J}\frac{\partial J}{\partial\tau}H_{0\alpha}-\frac{1}{c}\frac{\partial c}{\partial\tau}H_{0\alpha}-\frac{\partial g_{\alpha\beta}}{\partial\tau}H_0^{~\beta}+\frac{\sigma}{\varepsilon}H_{0\alpha}=0.
\]
Introducing the vector $\eta$ such that $\eta_3=0$ and
\begin{equation}\label{equationeta}
\frac{\partial \eta_\alpha}{\partial \tau}-\frac{1}{2}\frac{\partial g_{\alpha\beta}}{\partial\tau}\eta^\beta+c^{-1}\partial_3c\eta_\alpha=0,
\end{equation}
and denote
\[
H_{0\alpha} =A\varepsilon^{1/2}c\eta_\alpha.
\]
If we consider the equation \eqref{equationeta} as an ODE along a geodesic $\gamma$ with respect to the metric $g$, we can see that $\eta$ is a parallel vector field along $\gamma$ and orthogonal to $\gamma$. So one can take $c^{-2}|\eta|^2=1$.
We obtain the following equation for $A$,
\[
2\frac{\partial A}{\partial\tau}+\frac{1}{J}\frac{\partial J}{\partial\tau}A-\frac{1}{c}\frac{\partial c}{\partial\tau}A+\frac{\sigma}{\varepsilon}A=0.
\]
Solving the above equation along $\gamma$, we have
\[
A=CJ^{-1/2}c^{1/2}e^{-\int\frac{\sigma}{2\varepsilon}\mathrm{d}\tau}
\]
with some constant $C$ and
\[
I=e^{-\int\frac{\sigma}{2\varepsilon}\mathrm{d}\tau}.
\]
Since we need nontrivial solutions, we just take $C\equiv 1$.

So from the impedance map, knowing the value of $\varepsilon$ on $\Gamma$, we can recover the geodesic ray transform of $\frac{\sigma}{\varepsilon}$
\[
\int_\gamma \frac{\sigma}{\varepsilon}(\gamma(s))\mathrm{d}s,
\]
where $\gamma$ is an arbitrary geodesic connecting two  points $a,b\in\Gamma$. By the invertibility of local geodesic ray transform \cite{uhlmann2016inverse}, we can uniquely determine $\frac{\sigma}{\varepsilon}$ in $\kappa^{-1}([0,1])\cap\overline{\Omega}$.

\subsection{Recovery of $\varepsilon$}
Introducing the vector $\zeta$ such that
\[
\zeta^1=-cg^{-1/2}\eta_2,\quad \zeta^2=cg^{-1/2}\eta_1,\quad \zeta^3=0.
\]
So
\[
\begin{split}
c^{-2}|\zeta|^2=&c^{-2}(g_{11}\zeta^1\zeta^1+2g_{12}\zeta^1\zeta^2+g_{22}\zeta^2\zeta^2)\\
=&c^{-2}c^2g^{-1}J^{2}(g^{22}\eta_2\eta_2+2g^{12}\eta_1\eta_2+2g^{11}\eta_1\eta_1)\\
=&c^{-2}(g^{22}\eta_2\eta_2+2g^{12}\eta_1\eta_2+2g^{11}\eta_1\eta_1)=1.
\end{split}
\]
One can see that $\zeta$ is also parallel along $\gamma$.
Then \eqref{E0H0relation} implies
\[
E_0^{~\alpha}=-A\mu^{1/2}c\zeta^\alpha,\quad\alpha=1,2, 
\]
that is, $E_0=-A\mu^{1/2}c\zeta$. We conclude that the principal term in the amplitudes does not contain any information of $\varepsilon$ in the interior of the domain with $\varepsilon\mu,\frac{\sigma}{\varepsilon}$ already determined. In order to recover $\varepsilon$, we need to use the subprincipal term of the amplitudes.\\

%

Putting $m=1$ and $i=3$ in \eqref{GOequations}, we have
\[
\begin{split}
\mu H^{~3}_1&=e^{3\alpha\beta}E_{0\beta;\alpha}=-g^{-1/2}(E_{01;2}-E_{02;1}),\\
\varepsilon E_1^{~3}&=-e^{3\alpha\beta}H_{0\beta;\alpha}=g^{-1/2}(H_{01;2}-H_{02;1})+\sigma E_0^{~3}.
\end{split}
\]
Putting $m=2$ and $i=\alpha=1,2$, we obtain
\begin{align}
\mu H_2^{~1}&=-g^{-1/2}E_{22}-g^{-1/2}(E_{12;3}-E_{13;2}),\label{eqH21}\\
\mu H_2^{~2}&=g^{-1/2}E_{21}+g^{-1/2}(E_{11;3}-E_{13;1}),\label{eqH22}\\
\varepsilon E_2^{~1}&=g^{-1/2}H_{22}+g^{-1/2}(H_{12;3}-H_{13;2})+\sigma E_1^{~1},\label{eqE21}\\
\varepsilon E_2^{~2}&=-g^{-1/2}H_{21}-g^{-1/2}(H_{11;3}-H_{13;1})+\sigma E_1^{~2}.\label{eqE22}
\end{align}

Upon using $H_{1\alpha}=g_{\alpha\beta}H_{1}^{~\beta}$ in \eqref{eqH21} and \eqref{eqH22}, we have
\[
\begin{split}
\mu H_{21}=-g^{1/2}c^{-2} (E_2^{~2}+E_{1~;3}^{~2}-E_{13;}^{~\,~\,~2}),
\end{split}
\]
and similarly
\[
\mu H_{22}=g^{1/2}c^{-2} (E_2^{~1}+E_{1~;3}^{~1}-E_{13;}^{~\,~\,~1}).
\]
Inserting into \eqref{eqE21} and \eqref{eqE22}, we obtain
\begin{equation}\label{mainidentity1}
c^{-2}(E_{1~;3}^{~1}-E_{13;}^{~\,~\,1})+\mu g^{-1/2}(H_{12;3}-H_{13;2})+\mu\sigma E_1^{~1}=0,
\end{equation}
\begin{equation}\label{mainidentity2}
c^{-2}(E_{1~;3}^{~2}-E_{13;}^{~\,~\,2})-\mu g^{-1/2}(H_{11;3}-H_{13;1})+\mu\sigma E_1^{~2}=0.
\end{equation}
Taking \eqref{mainidentity1}$\times\eta^2$-\eqref{mainidentity2}$\times\eta^1$, we obtain
\begin{equation}\label{mainidentity}
\begin{split}
&\left(c^{-2}(E_{1~;3}^{~1}-E_{13;}^{~\,~\,1}+\mu g^{-1/2}(H_{12;3}-H_{13;2})+\mu\sigma E_1^{~1}\right)\eta^2\\
&\quad\quad\quad\quad\quad-\left(c^{-2}(E_{1~;3}^{~2}-E_{13;}^{~\,~\,2})-\mu g^{-1/2}(H_{11;3}-H_{13;1})+\mu\sigma E_1^{~2}\right)\eta^1=0.
\end{split}
\end{equation}
This would be the main equation to solve for the construction of the subprincipal terms $H_1,E_1$ of the amplitudes.\\

\noindent\textbf{Calculation of $E_{1~;3}^{~1}$ and $E_{1~;3}^{~2}$.}
We calculate, using \eqref{eqE21},
\begin{align*}
E_{1~;3}^{~1}=&\partial_3E_1^{~1}+\Gamma^1_{3\alpha}E_1^{~\alpha}\\
=&\partial_3\left(\frac{1}{\varepsilon}g^{-1/2}H_{12}+\frac{1}{\varepsilon}g^{-1/2}(H_{02;3}-H_{03;2})+\frac{\sigma}{\varepsilon}E_0^{~1}\right)+\Gamma_{31}^1\frac{1}{\varepsilon}g^{-1/2}(H_{12}+H_{02;3}-H_{03;2}+\frac{\sigma}{\varepsilon}E_0^{~1})\\
&-\Gamma_{32}^1\frac{1}{\varepsilon}g^{-1/2}(H_{11}+H_{01;3}-H_{03;1}-\frac{\sigma}{\varepsilon}E_0^{~2})+\Gamma_{33}^1\partial_\beta c\frac{1}{\varepsilon}g^{-1/2}(H_{01;2}-H_{02;1})\\
=&-\varepsilon^{-2}\partial_3\varepsilon g^{-1/2}(H_{12}+\partial_3H_{02})-\frac{1}{2}\varepsilon^{-1}g^{-3/2}\frac{\partial g}{\partial\tau}(H_{12}+\partial_3H_{02})\\
&+\varepsilon^{-1}g^{-1/2}\partial_3(H_{12}+\partial_3H_{02})+\partial_3(\sigma\varepsilon^{-1})E_0^{~1}+\frac{\sigma}{\varepsilon}\partial_3E_0^{~1}\\
&+\frac{1}{2}g^{1\gamma}\frac{\partial g_{1\gamma}}{\partial\tau}\left(\varepsilon^{-1}g^{-1/2}(H_{12}+\partial_3H_{02})+\frac{\sigma}{\varepsilon}E_0^{~1}\right)-\frac{1}{2}g^{1\gamma}\frac{\partial g_{2\gamma}}{\partial\tau}\left(\varepsilon^{-1}g^{-1/2}(H_{11}+\partial_3H_{01})-\frac{\sigma}{\varepsilon}E_0^{~2}\right)\\
&-cg^{1\beta}\partial_\beta c\frac{1}{\varepsilon}g^{-1/2}(H_{01;2}-H_{02;1})\\
=&\varepsilon^{-1}g^{-1/2}\left(-\varepsilon^{-1}\partial_3\varepsilon H_{12}-\frac{1}{2}g\frac{\partial g}{\partial\tau}H_{12}+\partial_3H_{12}+\frac{1}{2}g^{1\gamma}\frac{\partial g_{1\gamma}}{\partial\tau}H_{12}-\frac{1}{2}g^{1\gamma}\frac{\partial g_{2\gamma}}{\partial\tau}H_{11}\right)\\
&-g^{-1/2}\varepsilon^{-2}\partial_3\varepsilon\partial_3(A\varepsilon^{1/2}c\eta_2)-\frac{1}{2}\varepsilon^{-1}g^{-3/2}\frac{\partial g}{\partial\tau}\partial_3(A\varepsilon^{1/2}c\eta_2)+\varepsilon^{-1}g^{-1/2}\partial^2_3(A\varepsilon^{1/2}c\eta_2)\\&+\frac{1}{2}g^{1\gamma}\frac{\partial g_{1\gamma}}{\partial\tau}\varepsilon^{-1}g^{-1/2}\partial_3(A\varepsilon^{1/2}c\eta_2)-\frac{1}{2}g^{1\gamma}\frac{\partial g_{2\gamma}}{\partial\tau}\varepsilon^{-1}g^{-1/2}\partial_3(A\varepsilon^{1/2}c\eta_1)+\frac{\sigma}{\varepsilon}\partial_3(Ac\varepsilon^{-1/2}g^{-1/2}\eta_2)\\
&-cg^{1\beta}\partial_\beta c\varepsilon^{-1}\left(\partial_2(A\varepsilon^{1/2}c\eta_1)-\partial_1(A\varepsilon^{1/2}c\eta_2)\right)+\varepsilon^{-1/2}R,
\end{align*}
where we have used \eqref{eqE1} and \eqref{eqE2}. Here and in the following $R$, which may vary from step to step, denotes terms that do not depend on $\varepsilon$, with $\varepsilon\mu$ and $\frac{\sigma}{\varepsilon}$ already known.

We will use the identities
\[
\begin{split}
\partial_3(A\varepsilon^{1/2}c\eta_2)=&\partial_3(J^{-1/2}\varepsilon^{1/2}c^{3/2}I\eta_2)\\
=&-\frac{1}{2}J^{-3/2}\partial_3J\varepsilon^{1/2}c^{3/2}I\eta_2+J^{-1/2}\varepsilon^{1/2}\partial_3\log\sqrt{\varepsilon}c^{3/2}I\eta_2\\
&+\frac{3}{2}J^{-1/2}\varepsilon^{1/2}c^{1/2}\partial_3c I\eta_2-\frac{1}{2}J^{-1/2}\varepsilon^{1/2}c^{3/2}\frac{\sigma}{\varepsilon}I\eta_2+J^{-1/2}\varepsilon^{1/2}c^{3/2}I\partial_3\eta_2,
\end{split}
\]
\[
\begin{split}
\partial^2_{33}(A\varepsilon^{1/2}c\eta_2)=-&J^{-3/2}\partial_3J\varepsilon^{1/2}\partial_3\log\sqrt{\varepsilon}c^{3/2}I\eta_2+J^{-1/2}\varepsilon^{1/2}(\partial_3\log\sqrt{\varepsilon})^2c^{3/2}I\eta_2\\
&+J^{-1/2}\varepsilon^{1/2}\partial_{33}^2\log\sqrt{\varepsilon}c^{3/2}I\eta_2+3J^{-1/2}\varepsilon^{1/2}\partial_3\log\sqrt{\varepsilon}c^{1/2}\partial_3cI\eta_2\\
&+2J^{-1/2}\varepsilon^{1/2}c^{3/2}\partial_3\log\sqrt{\varepsilon}I\partial_3\eta_2-\frac{\sigma}{\varepsilon}J^{-1/2}\varepsilon^{1/2}c^{3/2}\partial_3\log\sqrt{\varepsilon}I\eta_2+\varepsilon^{1/2}R.
\end{split}
\]

Now we have
\begin{align}\label{eqE113}
E_{1~;3}^{~1}\notag
=&-\varepsilon^{-1}g^{-1/2}\left(\varepsilon^{-1}\partial_3\varepsilon H_{12}+\frac{1}{2}g^{-1}\frac{\partial g}{\partial\tau}H_{12}-\partial_3H_{12}-\frac{1}{2}g^{1\gamma}\frac{\partial g_{1\gamma}}{\partial\tau}H_{12}+\frac{1}{2}g^{1\gamma}\frac{\partial g_{2\gamma}}{\partial\tau}H_{11}\right)\notag\\
&-g^{-1/2}(\partial_3\log\sqrt{\varepsilon})^2J^{-1/2}c^{3/2}I\eta_2\\
&-\frac{1}{2}\varepsilon^{-1/2}g^{-3/2}g^{\beta\gamma}\frac{\partial g_{\beta\gamma}}{\partial\tau}J^{-1/2}\partial_3\log\sqrt{\varepsilon}c^{3/2}I\eta_2-g^{-1/2}\varepsilon^{-1/2}\partial_3\log\sqrt{\varepsilon}J^{-1/2}c^{1/2}\partial_3cI\eta_2\notag\\
&+g^{-1/2}J^{-1/2}\varepsilon^{-1/2}\partial_{33}^2\log\sqrt{\varepsilon}c^{3/2}I\eta_2\notag\\
&+\frac{1}{2}g^{1\gamma}\frac{\partial g_{1\gamma}}{\partial\tau}g^{-1/2}\varepsilon^{-1/2}\partial_3\log\sqrt{\varepsilon}J^{-1/2}c^{3/2}I\eta_2-\frac{1}{2}g^{1\gamma}\frac{\partial g_{2\gamma}}{\partial\tau}g^{-1/2}\varepsilon^{-1/2}\partial_3\log\sqrt{\varepsilon}J^{-1/2}c^{3/2}I\eta_1\notag\\
&-\frac{\sigma}{\varepsilon}\varepsilon^{-1/2}g^{-1/2}J^{-1/2}c^{3/2}\partial_3\log\sqrt{\varepsilon}I\eta_2\notag\\
&-\varepsilon^{-1/2}g^{-1/2}J^{-1/2}c^{5/2}g^{1\beta}\partial_\beta c(\partial_2\log\sqrt{\varepsilon}I\eta_1-\partial_1\log\sqrt{\varepsilon}I\eta_2)+\varepsilon^{-1/2}R.\notag
\end{align}
Similar calculation leads to
\begin{align}\label{eqE123}
-E_{1~;3}^{~2}\notag
=&\varepsilon^{-1}g^{-1/2}\left(-\varepsilon^{-1}\partial_3\varepsilon H_{11}-\frac{1}{2}g^{-1}\frac{\partial g}{\partial\tau}H_{11}+\partial_3H_{11}-\frac{1}{2}g^{2\gamma}\frac{\partial g_{1\gamma}}{\partial\tau}H_{12}+\frac{1}{2}g^{2\gamma}\frac{\partial g_{2\gamma}}{\partial\tau}H_{11}\right)\notag\\
&-g^{-1/2}\varepsilon^{-1/2}(\partial_3\log\sqrt{\varepsilon})^2J^{-1/2}c^{3/2}I\eta_1\\
&-\frac{1}{2}\varepsilon^{-1/2}g^{-3/2}g^{\beta\gamma}\frac{\partial g_{\beta\gamma}}{\partial\tau}J^{-1/2}\partial_3\log\sqrt{\varepsilon}c^{3/2}I\eta_1-g^{-1/2}\varepsilon^{-1/2}\partial_3\log\sqrt{\varepsilon}J^{-1/2}c^{1/2}\partial_3cI\eta_1\notag\\
&+g^{-1/2}J^{-1/2}\varepsilon^{-1/2}\partial_{33}^2\log\sqrt{\varepsilon}c^{3/2}I\eta_1\notag\\
&-\frac{1}{2}g^{2\gamma}\frac{\partial g_{1\gamma}}{\partial\tau}g^{-1/2}\varepsilon^{-1/2}\partial_3\log\sqrt{\varepsilon}J^{-1/2}c^{3/2}I\eta_2+\frac{1}{2}g^{2\gamma}\frac{\partial g_{2\gamma}}{\partial\tau}g^{-1/2}\varepsilon^{-1/2}\partial_3\log\sqrt{\varepsilon}J^{-1/2}c^{3/2}I\eta_1\notag\\
&-\frac{\sigma}{\varepsilon}\varepsilon^{-1/2}g^{-1/2}J^{-1/2}c^{3/2}\partial_3\log\sqrt{\varepsilon}I\eta_1\notag\\
&-\varepsilon^{-1/2}g^{-1/2}J^{-1/2}c^{5/2}g^{2\beta}\partial_\beta c(\partial_1\log\sqrt{\varepsilon}I\eta_2-\partial_2\log\sqrt{\varepsilon}I\eta_1)+\varepsilon^{-1/2}R.\notag
\end{align}
Here we have used the identity
\[
\partial g=2c\partial cJ^2+2c^2J\partial J.
\]

\noindent\textbf{Calculation of $\mu\sigma E_1^{~1}\eta^2-\mu\sigma E_1^{~2}\eta^1$.}
We calculate
\[
\begin{split}
&\mu\sigma E_1^{~1}\eta^2-\mu\sigma E_1^{~2}\eta^1\\
=&\mu\sigma \frac{1}{\varepsilon}g^{-1/2}\partial_3H_{02}\eta^2+\mu\sigma\frac{1}{\epsilon}g^{-1/2}\partial_3H_{01}\eta^1+\varepsilon^{-1/2}R\\
=&\mu\sigma \frac{1}{\varepsilon}g^{-1/2}\varepsilon^{1/2}J^{-1/2}c^{3/2}\partial_3\log\sqrt{\varepsilon}(\eta_2\eta^2+\eta_1\eta^1)+\varepsilon^{-1/2}R\\
=&c^{-2}\frac{\sigma}{\varepsilon}g^{-1/2}\varepsilon^{-1/2}J^{-1/2}c^{3/2}\partial_3\log\sqrt{\varepsilon}(\eta_2\eta^2+\eta_1\eta^1)+\varepsilon^{-1/2}R.
\end{split}
\]

\noindent\textbf{Calculation of $E_{13;}^{\phantom~\quad~1}$ and $E_{13;}^{\phantom~\quad~2}$.}
Putting $m=1$, $i=3$ in \eqref{GOequations}, we have
\[
\begin{split}
E_{13}=&\frac{1}{\varepsilon}g^{-1/2}g_{33}(H_{01;2}-H_{02;1})\\
=&\frac{1}{\varepsilon}g^{-1/2}c^2(\partial_2H_{01}-\Gamma^\gamma_{12}H_{0\gamma}-\partial_1H_{02}+\Gamma^\gamma_{12}H_{0\gamma})\\
=&\frac{1}{\varepsilon}g^{-1/2}c^2\Bigg(\partial_2A\varepsilon^{1/2}c\eta_1+A\partial_2\varepsilon^{1/2}c\eta_1+A\varepsilon^{1/2}\partial_2c\eta_1+A\varepsilon^{1/2}c\partial_2\eta_1\\
&\quad\quad\quad-\partial_1A\varepsilon^{1/2}c\eta_2-A\partial_1\varepsilon^{1/2}c\eta_2-A\varepsilon^{1/2}\partial_1c\eta_2-A\varepsilon^{1/2}c\partial_1\eta_2
\Bigg)\\
=&\frac{1}{\sqrt{\varepsilon}}g^{-1/2}c^2\Bigg(\partial_2Ac\eta_1+A\partial_2\log\sqrt{\varepsilon}c\eta_1+A\partial_2c\eta_1+Ac\partial_2\eta_1\\
&\quad\quad\quad-\partial_1Ac\eta_2-A\partial_1\log\sqrt{\varepsilon}c\eta_2-A\partial_1c\eta_2-Ac\partial_1\eta_2
\Bigg).
\end{split}
\]
We calculate
\begin{align*}
g^{\alpha\beta}E_{13;\beta}&=g^{\alpha\beta}\partial_\beta E_{13}-g^{\alpha\beta}\Gamma_{3\beta}^\gamma E_{1\gamma}-g^{\alpha\beta}\Gamma^3_{3\beta}E_{13}\\
&=-\frac{1}{\sqrt{\varepsilon}}g^{\alpha\beta}g^{-1/2}c^2\partial_\beta\log\sqrt{\varepsilon}\Bigg(\partial_2Ac\eta_1+A\partial_2\log\sqrt{\varepsilon}c\eta_1+A\partial_2c\eta_1+Ac\partial_2\eta_1\\
&\quad\quad\quad-\partial_1Ac\eta_2-A\partial_1\log\sqrt{\varepsilon}c\eta_2-A\partial_1c\eta_2-Ac\partial_1\eta_2
\Bigg)\\
&+\frac{1}{\sqrt{\varepsilon}}g^{\alpha\beta}\partial_\beta(g^{-1/2}Ac^3\eta_1)\partial_2\log\sqrt{\varepsilon}-\frac{1}{\sqrt{\varepsilon}}g^{\alpha\beta}\partial_\beta(g^{-1/2}Ac^3\eta_2)\partial_1\log\sqrt{\varepsilon}\\
&+\frac{1}{\sqrt{\varepsilon}}g^{\alpha\beta}(g^{-1/2}Ac^3\eta_1)\partial_\beta\partial_2\log\sqrt{\varepsilon}-\frac{1}{\sqrt{\varepsilon}}g^{\alpha\beta}(g^{-1/2}Ac^3\eta_2)\partial_\beta\partial_1\log\sqrt{\varepsilon}\\
&-g^{\alpha\beta}c^{-1}\partial_\beta c\frac{1}{\sqrt{\varepsilon}}c^2g^{-1/2}(Ac\partial_2\log\sqrt{\varepsilon}\eta_1-Ac\partial_1\log\sqrt{\varepsilon}\eta_2)\\
&-\frac{1}{2}g^{\alpha\beta}\frac{\partial g_{\beta\delta}}{\partial\tau}E_1^{~\delta}
\end{align*}
Note that
\[
\begin{split}
\frac{1}{2}g^{\alpha\beta}\frac{\partial g_{\beta\delta}}{\partial\tau}E_1^{~\delta}=&\frac{1}{2}g^{\alpha\beta}\frac{\partial g_{\beta1}}{\partial\tau}E_1^{~1}+\frac{1}{2}g^{\alpha\beta}\frac{\partial g_{\beta2}}{\partial\tau}E_1^{~2}\\
=&\frac{1}{2}g^{\alpha\beta}\frac{\partial g_{\beta1}}{\partial\tau}\left(\frac{1}{\varepsilon}g^{-1/2}(H_{12}+H_{02;3}-H_{03;2})+\frac{\sigma}{\varepsilon}E_0^{~1}\right)\\
&-\frac{1}{2}g^{\alpha\beta}\frac{\partial g_{\beta2}}{\partial\tau}\frac{1}{\varepsilon}\left(g^{-1/2}(H_{11}+H_{01;3}-H_{03;1})-\frac{\sigma}{\varepsilon}E_0^{~2}\right)\\
=&\frac{1}{2}g^{\alpha\beta}\frac{\partial g_{\beta1}}{\partial\tau}\frac{1}{\varepsilon}g^{-1/2}H_{12}-\frac{1}{2}g^{\alpha\beta}\frac{\partial g_{\beta2}}{\partial\tau}\frac{1}{\varepsilon}g^{-1/2}H_{11}\\
&+\frac{1}{2}g^{\alpha\beta}\frac{\partial g_{\beta1}}{\partial\tau}\varepsilon^{-1/2}g^{-1/2}c^{3/2}J^{-1/2}\partial_3\log\sqrt{\varepsilon}I\eta_2\\
&-\frac{1}{2}g^{\alpha\beta}\frac{\partial g_{\beta2}}{\partial\tau}\varepsilon^{-1/2}g^{-1/2}c^{3/2}J^{-1/2}\partial_3\log\sqrt{\varepsilon}I\eta_1+\varepsilon^{-1/2}R.
\end{split}
\]
So
\begin{align*}
g^{\alpha\beta}E_{13;\beta}&=-g^{\alpha\beta}g^{-1/2}\varepsilon^{-1/2}c^{7/2}J^{-1/2}\partial_\beta\log\sqrt{\varepsilon}\Bigg(-\frac{1}{2}J^{-1}\partial_2JI\eta_1+\frac{3}{2}c^{-1}\partial_2cI\eta_1+\partial_2\log\sqrt{\varepsilon}I\eta_1+\partial_2(I\eta_1)\\
&\quad\quad+\frac{1}{2}J^{-1}\partial_1JI\eta_2-\frac{3}{2}c^{-1}\partial_1cI\eta_2-\partial_1\log\sqrt{\varepsilon}I\eta_2-\partial_1(I\eta_2)\Bigg)\\
&+g^{\alpha\beta}g^{-1/2}\varepsilon^{-1/2}c^{7/2}J^{-1/2}\partial_2\log\sqrt{\varepsilon}\left(\frac{5}{2}c^{-1}\partial_\beta cI\eta_1-\frac{3}{2}J^{-1}\partial_\beta JI\eta_1+\partial_\beta(I\eta_1)\right)\\
&-g^{\alpha\beta}g^{-1/2}\varepsilon^{-1/2}c^{7/2}J^{-1/2}\partial_1\log\sqrt{\varepsilon}\left(\frac{5}{2}c^{-1}\partial_\beta cI\eta_2-\frac{3}{2}J^{-1}\partial_\beta JI\eta_2+\partial_\beta(I\eta_2)\right)\\
&-g^{\alpha\beta}g^{-1/2}\varepsilon^{-1/2}c^{7/2}J^{-1/2}(-\partial_\beta\partial_2\log\sqrt{\varepsilon}I\eta_1+\partial_\beta\partial_1\log\sqrt{\varepsilon}I\eta_2)\\
&-g^{\alpha\beta}g^{-1/2}\varepsilon^{-1/2}c^{7/2}J^{-1/2}c^{-1}\partial_\beta c(\partial_2\log\sqrt{\varepsilon}I\eta_1-\partial_1\log\sqrt{\varepsilon}I\eta_2)\\
&-\frac{1}{2}g^{\alpha\beta}\frac{\partial g_{\beta1}}{\partial\tau}\frac{1}{\varepsilon}g^{-1/2}H_{12}+\frac{1}{2}g^{\alpha\beta}\frac{\partial g_{\beta2}}{\partial\tau}\frac{1}{\varepsilon}g^{-1/2}H_{11}\\
&-\frac{1}{2}g^{\alpha\beta}\frac{\partial g_{\beta1}}{\partial\tau}\varepsilon^{-1/2}g^{-1/2}c^{3/2}J^{-1/2}\partial_3\log\sqrt{\varepsilon}I\eta_2\\
&+\frac{1}{2}g^{\alpha\beta}\frac{\partial g_{\beta2}}{\partial\tau}\varepsilon^{-1/2}g^{-1/2}c^{3/2}J^{-1/2}\partial_3\log\sqrt{\varepsilon}I\eta_1+\varepsilon^{-1/2}R.
\end{align*}

Taking $\alpha=1$ and then $\alpha=2$, we have
\begin{align*}
&-c^{-2}g^{1\beta}E_{13;\beta}\\
&=g^{-1/2}\varepsilon^{-1/2}c^{3/2}J^{-1/2}\Bigg(-\frac{1}{2}g^{11}\partial_1\log\sqrt{\varepsilon}\partial_2JI\eta_1-g^{11}\partial_1\log\sqrt{\varepsilon}\partial_1JI\eta_2+g^{12}\partial_2\log\sqrt{\varepsilon}\partial_2JI\eta_1\\
&\quad\quad\quad\quad\quad+\frac{1}{2}g^{12}\partial_2\log\sqrt{\varepsilon}\partial_1JI\eta_2+\frac{3}{2}g^{11}\partial_2\log\sqrt{\varepsilon}\partial_1JI\eta_1-\frac{3}{2}g^{12}\partial_1\log\sqrt{\varepsilon}\partial_2J\eta_2\Bigg)\\
&+g^{-1/2}\varepsilon^{-1/2}c^{3/2}J^{-1/2}\left(-g^{11}\partial_1\partial_2\log\sqrt{\varepsilon}I\eta_1+g^{11}\partial_1\partial_1\log\sqrt{\varepsilon}I\eta_2-g^{12}\partial_2\partial_2\log\sqrt{\varepsilon}I\eta_1+g^{12}\partial_2\partial_1\log\sqrt{\varepsilon}I\eta_2\right)\\
&+g^{-1/2}\varepsilon^{-1/2}c^{3/2}J^{-1/2}\Big(g^{11}\partial_1\log\sqrt{\varepsilon}\partial_2\log\sqrt{\varepsilon}I\eta_1-g^{11}\partial_1\log\sqrt{\varepsilon}\partial_1\log\sqrt{\varepsilon}I\eta_2\\
&\quad\quad\quad\quad\quad\quad\quad\quad\quad+g^{12}\partial_2\log\sqrt{\varepsilon}\partial_2\log\sqrt{\varepsilon}I\eta_1-g^{12}\partial_2\log\sqrt{\varepsilon}\partial_1\log\sqrt{\varepsilon}I\eta_2\Big)\\
&+g^{-1/2}\varepsilon^{-1/2}c^{3/2}J^{-1/2}\left(\frac{3}{2}g^{11}\partial_2c\partial_1\log\sqrt{\varepsilon}I\eta_1-\frac{3}{2}g^{11}\partial_2\log\sqrt{\varepsilon}\partial_1cI\eta_1\right.\\
&\quad\quad\quad\quad\quad\quad\quad\quad\quad\left.-\frac{3}{2}g^{12}\partial_2\log\sqrt{\varepsilon}\partial_1cI\eta_2+\frac{3}{2}g^{12}\partial_1\log\sqrt{\varepsilon}\partial_2cI\eta_2\right)\\
&+g^{-1/2}\varepsilon^{-1/2}c^{3/2}J^{-1/2}(g^{11}\partial_1\log\sqrt{\varepsilon}\partial_2(I\eta_1)-g^{12}\partial_2\log\sqrt{\varepsilon}\partial_1(I\eta_2)-g^{11}\partial_2\log\sqrt{\varepsilon}\partial_1(I\eta_1)+g^{12}\partial_1\log\sqrt{\varepsilon}\partial_2(I\eta_2))\\
&+\frac{1}{2}g^{1\beta}\frac{\partial g_{\beta1}}{\partial\tau}\frac{1}{\varepsilon}g^{-1/2}H_{12}-\frac{1}{2}g^{1\beta}\frac{\partial g_{\beta2}}{\partial\tau}\frac{1}{\varepsilon}g^{-1/2}H_{11}\\
&+\frac{1}{2}g^{1\beta}\frac{\partial g_{\beta1}}{\partial\tau}\varepsilon^{-1/2}g^{-1/2}c^{3/2}J^{-1/2}\partial_3\log\sqrt{\varepsilon}I\eta_2-\frac{1}{2}g^{1\beta}\frac{\partial g_{\beta2}}{\partial\tau}\varepsilon^{-1/2}g^{-1/2}c^{3/2}J^{-1/2}\partial_3\log\sqrt{\varepsilon}I\eta_1+\varepsilon^{-1/2}R,
\end{align*}
and
\begin{align*}
&c^{-2}g^{2\beta}E_{13;\beta}\\
&=g^{-1/2}\varepsilon^{-1/2}c^{3/2}J^{-1/2}\Bigg(-\frac{1}{2}g^{22}\partial_2\log\sqrt{\varepsilon}\partial_1JI\eta_2-g^{22}\partial_2\log\sqrt{\varepsilon}\partial_2JI\eta_1+g^{12}\partial_1\log\sqrt{\varepsilon}\partial_1JI\eta_2\\
&\quad\quad\quad\quad\quad+\frac{1}{2}g^{12}\partial_1\log\sqrt{\varepsilon}\partial_2JI\eta_1+\frac{3}{2}g^{22}\partial_1\log\sqrt{\varepsilon}\partial_2JI\eta_2-\frac{3}{2}g^{12}\partial_2\log\sqrt{\varepsilon}\partial_1J\eta_1\Bigg)\\
&+g^{-1/2}\varepsilon^{-1/2}c^{3/2}J^{-1/2}\left(-g^{22}\partial_2\partial_1\log\sqrt{\varepsilon}I\eta_2+g^{22}\partial_2\partial_2\log\sqrt{\varepsilon}I\eta_1-g^{21}\partial_1\partial_1\log\sqrt{\varepsilon}I\eta_2+g^{21}\partial_1\partial_2\log\sqrt{\varepsilon}I\eta_1\right)\\
&+g^{-1/2}\varepsilon^{-1/2}c^{3/2}J^{-1/2}\Big(g^{22}\partial_2\log\sqrt{\varepsilon}\partial_1\log\sqrt{\varepsilon}I\eta_2-g^{22}\partial_2\log\sqrt{\varepsilon}\partial_2\log\sqrt{\varepsilon}I\eta_1\\
&\quad\quad\quad\quad\quad\quad\quad\quad\quad+g^{21}\partial_1\log\sqrt{\varepsilon}\partial_1\log\sqrt{\varepsilon}I\eta_2-g^{21}\partial_1\log\sqrt{\varepsilon}\partial_2\log\sqrt{\varepsilon}I\eta_1\Big)\\
&+g^{-1/2}\varepsilon^{-1/2}c^{3/2}J^{-1/2}\left(-\frac{3}{2}g^{21}\partial_2c\partial_1\log\sqrt{\varepsilon}I\eta_1+\frac{3}{2}g^{21}\partial_2\log\sqrt{\varepsilon}\partial_1cI\eta_1\right.\\
&\quad\quad\quad\quad\quad\quad\quad\quad\quad\left.+\frac{3}{2}g^{22}\partial_2\log\sqrt{\varepsilon}\partial_1cI\eta_2-\frac{3}{2}g^{22}\partial_1\log\sqrt{\varepsilon}\partial_2cI\eta_2\right)\\
&+g^{-1/2}\varepsilon^{-1/2}c^{3/2}J^{-1/2}(g^{22}\partial_2\log\sqrt{\varepsilon}\partial_1(I\eta_2)-g^{21}\partial_1\log\sqrt{\varepsilon}\partial_2(I\eta_1)-g^{22}\partial_1\log\sqrt{\varepsilon}\partial_2(I\eta_2)+g^{21}\partial_2\log\sqrt{\varepsilon}\partial_1(I\eta_1))\\
&+\frac{1}{2}g^{2\beta}\frac{\partial g_{\beta 2}}{\partial\tau}\frac{1}{\varepsilon}g^{-1/2}H_{11}-\frac{1}{2}g^{2\beta}\frac{\partial g_{\beta1}}{\partial\tau}\frac{1}{\varepsilon}g^{-1/2}H_{12}\\
&+\frac{1}{2}g^{2\beta}\frac{\partial g_{\beta2}}{\partial\tau}\varepsilon^{-1/2}g^{-1/2}c^{3/2}J^{-1/2}\partial_3\log\sqrt{\varepsilon}I\eta_1-\frac{1}{2}g^{2\beta}\frac{\partial g_{\beta1}}{\partial\tau}\varepsilon^{-1/2}g^{-1/2}c^{3/2}J^{-1/2}\partial_3\log\sqrt{\varepsilon}I\eta_2+\varepsilon^{-1/2}R.
\end{align*}

Noticing
\[
\begin{split}
&g^{1\beta}\frac{\partial g_{1\beta}}{\partial\tau}\partial_3\log\sqrt{\varepsilon}\eta_2\eta^2-g^{1\beta}\frac{\partial g_{2\beta}}{\partial\tau}\partial_3\log\sqrt{\varepsilon}\eta_1\eta^2-g^{2\beta}\frac{\partial g_{1\beta}}{\partial\tau}\partial_3\log\sqrt{\varepsilon}\eta_2\eta^1+g^{2\beta}\frac{\partial g_{2\beta}}{\partial\tau}\partial_3\log\sqrt{\varepsilon}\eta_1\eta^1\\
=&(g_{\alpha\beta}\frac{\partial g_{\alpha\beta}}{\partial\tau}|\eta|^2-\frac{\partial g_{\alpha\beta}}{\partial\tau}\eta^\alpha\eta^\beta)\partial_3\log\sqrt{\varepsilon},
\end{split}
\]
we obtain
\begin{equation}\label{eqE1312combined}
\begin{split}
&g^{1/2}\varepsilon^{1/2}c^{-3/2}J^{1/2}(-c^{-2}E_{13;}^{~\,~\,~1}\eta^2+c^{-2}E_{13;}^{~\,~\,~2}\eta^1)\\
=&g^{1/2}\varepsilon^{1/2}c^{-3/2}J^{1/2}(c^{-2}g^{1\beta}E_{13;\beta}\eta^2-c^{-2}g^{2\beta}E_{13;\beta}\eta^1)\\
=&\frac{1}{2}c^{-3/2}J^{1/2}\varepsilon^{-1/2}g^{1\beta}\frac{\partial g_{\beta1}}{\partial\tau}H_{12}\eta^2-\frac{1}{2}c^{-3/2}J^{1/2}\varepsilon^{-1/2}g^{1\beta}\frac{\partial g_{\beta2}}{\partial\tau}H_{11}\eta^2\\
&-\frac{1}{2}c^{-3/2}J^{1/2}\varepsilon^{-1/2}g^{2\beta}\frac{\partial g_{\beta1}}{\partial\tau}H_{12}\eta^1+\frac{1}{2}c^{-3/2}J^{1/2}\varepsilon^{-1/2}g^{2\beta}\frac{\partial g_{\beta2}}{\partial\tau}H_{11}\eta^1\\
&-\partial_\alpha\partial_\beta\log\sqrt{\varepsilon}I\eta^\alpha\eta^\beta+g^{\alpha\beta}\partial_\alpha\partial_\beta\log\sqrt{\varepsilon}I|\eta|^2\\
&+\partial_\alpha\log\sqrt{\varepsilon}\partial_\beta\log\sqrt{\varepsilon}I\eta^\alpha\eta^\beta-g^{\alpha\beta}\partial_\alpha \log\sqrt{\varepsilon}\partial_\beta\log\sqrt{\varepsilon}I|\eta|^2\\
&+J^{-1}\partial_\alpha J\partial_\beta\log\sqrt{\varepsilon}\eta^\alpha\eta^\beta-g^{\alpha\beta}J^{-1}\partial_\alpha J\partial_\beta\log\sqrt{\varepsilon}|\eta|^2\\
&+\left(-g^{11}\partial_2\log\sqrt{\varepsilon}\partial_1(I\eta_1)+g^{12}\partial_1\log\sqrt{\varepsilon}\partial_2(I\eta_2)+g^{11}\partial_1\log\sqrt{\varepsilon}\partial_2(I\eta_1)-g^{12}\partial_2\log\sqrt{\varepsilon}\partial_1(I\eta_2)\right)\eta^2\\
&+\left(-g^{22}\partial_1\log\sqrt{\varepsilon}\partial_2(I\eta_2)+g^{12}\partial_2\log\sqrt{\varepsilon}\partial_1(I\eta_1)+g^{22}\partial_2\log\sqrt{\varepsilon}\partial_1(I\eta_2)-g^{12}\partial_1\log\sqrt{\varepsilon}\partial_2(I\eta_1)\right)\eta^1\\
&+\frac{1}{2}(g_{\alpha\beta}\frac{\partial g_{\alpha\beta}}{\partial\tau}|\eta|^2-\frac{\partial g_{\alpha\beta}}{\partial\tau}\eta^\alpha\eta^\beta)\partial_3\log\sqrt{\varepsilon}+\varepsilon^{-1/2}R.
\end{split}
\end{equation}

\noindent\textbf{Calculation of $(H_{12;3}-H_{13;2})$ and $(H_{11;3}-H_{13;1})$.}
We calculate
\begin{equation}\label{equation22}
\begin{split}
\mu g^{-1/2}(H_{12;3}-H_{13;2})=&\mu g^{-1/2}(\partial_3H_{12}-\partial_2H_{13})\\
=&\mu g^{-1/2}(\partial_3H_{12}+\partial_2(c^2\frac{1}{\mu}g^{-1/2}\partial_2E_{01}-c^2\frac{1}{\mu}g^{-1/2}\partial_1E_{02}))\\
=&\mu g^{-1/2}\partial_3H_{12}+g^{-1/2}\partial_2(c^2g^{-1/2}\partial_2E_{01})-\mu^{-1}\partial_2\mu g^{-1/2}c^2g^{-1/2}\partial_2E_{01}\\
&\quad\quad-g^{-1/2}\partial_2(c^2g^{-1/2}\partial_1E_{02})+\mu^{-1}\partial_2\mu g^{-1/2}c^2g^{-1/2}\partial_1E_{02},
\end{split}
\end{equation}
\[
\begin{split}
\mu g^{-1/2}(H_{11;3}-H_{13;1})=&\mu g^{-1/2}(\partial_3H_{11}-\partial_1H_{13})\\
=&\mu g^{-1/2}(\partial_3H_{11}+\partial_1(c^2\frac{1}{\mu}g^{-1/2}\partial_2E_{01}-c^2\frac{1}{\mu}g^{-1/2}\partial_1E_{02}))\\
=&\mu g^{-1/2}\partial_3H_{11}+g^{-1/2}\partial_1(c^2g^{-1/2}\partial_2E_{01})-\mu^{-1}\partial_1\mu g^{-1/2}c^2g^{-1/2}\partial_2E_{01}\\
&\quad\quad-g^{-1/2}\partial_1(c^2g^{-1/2}\partial_1E_{02})+\mu^{-1}\partial_1\mu g^{-1/2}c^2g^{-1/2}\partial_1E_{02}.
\end{split}
\]
Note
\[
\begin{split}
E_{01}=&-A\mu^{1/2}c\zeta_1=-A\mu^{1/2}c(g_{11}\zeta^1+g_{12}\zeta^2)\\
=&-A\mu^{1/2}c(-cg_{11}g^{-1/2}\eta_2+cg_{12}g^{-1/2}\eta_1)\\
=&A\mu^{1/2}c^2J^2g^{-1/2}(g^{22}\eta_2+g^{12}\eta_1)\\
=&\varepsilon^{-1/2}c^{3/2}J^{3/2}g^{-1/2}I(g^{22}\eta_2+g^{12}\eta_1).
\end{split}
\]
So we have
\[
\begin{split}
&c^2g^{-1/2}\partial_2E_{01}\\
=&c^2g^{-1/2}\partial_2(J^{-3/2}c^{3/2}\varepsilon^{-1/2}g^{-1/2}(g^{22}I\eta_2+g^{12}I\eta_1))\\
=&\varepsilon^{-1/2}c^{3/2}J^{-1/2}(-\partial_2\log\sqrt{\varepsilon}+\frac{3}{2}c^{-1}\partial_2c+\frac{3}{2}J^{-1}\partial_2J-\frac{1}{2}g^{-1}\partial_2g)(g^{22}I\eta_2+g^{12}I\eta_1)\\
&+\varepsilon^{-1/2}c^{3/2}J^{-1/2}(\partial_2g^{22}I\eta_2+g^{22}\partial_2(I\eta_2)+\partial_2g^{12}I\eta_1+g^{12}\partial_2(I\eta_1)).
\end{split}
\]
Similarly, we also have
\[
\begin{split}
&c^2g^{-1/2}\partial_1E_{02}\\
=&-c^2g^{-1/2}\partial_1(J^{-3/2}c^{3/2}\varepsilon^{-1/2}g^{-1/2}(g^{12}I\eta_2+g^{11}I\eta_1))\\
=&-\varepsilon^{-1/2}c^{3/2}J^{-1/2}(-\partial_1\log\sqrt{\varepsilon}+\frac{3}{2}c^{-1}\partial_1c+\frac{3}{2}J^{-1}\partial_1J-\frac{1}{2}g^{-1}\partial_1g)(g^{12}I\eta_2+g^{11}I\eta_1)\\
&-\varepsilon^{-1/2}c^{3/2}J^{-1/2}(\partial_1g^{12}I\eta_2+g^{12}\partial_1(I\eta_2)+\partial_1g^{11}I\eta_1+g^{11}\partial_2(I\eta_1)).
\end{split}
\]

We proceed with the calculation
\begin{equation}\label{equation23}
\begin{split}
&g^{-1/2}\partial_2(c^2g^{-1/2}\partial_2E_{01})\\
=&-\varepsilon^{-1/2}c^{3/2}J^{-1/2}g^{-1/2}\partial_2\log\sqrt{\varepsilon}(-\partial_2\log\sqrt{\varepsilon}+\frac{3}{2}c^{-1}\partial_2c+\frac{3}{2}J^{-1}\partial_2J-\frac{1}{2}g^{-1}\partial_2g)(g^{22}I\eta_2+g^{12}I\eta_1)\\
&-\frac{3}{2}\varepsilon^{-1/2}c^{3/2}J^{-1/2}g^{-1/2}c^{-1}\partial_2c\partial_2\log\sqrt{\varepsilon}(g^{22}I\eta_2+g^{12}I\eta_1)\\
&+\frac{1}{2}\varepsilon^{-1/2}c^{3/2}J^{-1/2}g^{-1/2}J^{-1}\partial_2J\partial_2\log\sqrt{\varepsilon}(g^{22}I\eta_2+g^{12}I\eta_1)\\
&-\varepsilon^{-1/2}c^{3/2}J^{-1/2}g^{-1/2}\partial_2\partial_2\log\sqrt{\varepsilon}I\eta^2\\
&-2\varepsilon^{-1/2}c^{3/2}J^{-1/2}g^{-1/2}\partial_2\log\sqrt{\varepsilon}(\partial_2g^{22}I\eta_2+g^{22}\partial_2(I\eta_2)+\partial_2g^{12}I\eta_1+g^{12}\partial_2(I\eta_1))+\varepsilon^{-1/2}R.
\end{split}
\end{equation}

\begin{equation}\label{equation24}
\begin{split}
&-\mu^{-1}\partial_2\mu g^{-1/2}c^2g^{-1/2}\partial_2E_{01}\\
=&(2\partial_2\log\sqrt{\varepsilon}+2c^{-1}\partial_2c)g^{-1/2}c^2g^{-1/2}\partial_2E_{01}\\
=&2\partial_2\log\sqrt{\varepsilon}g^{-1/2}\varepsilon^{-1/2}c^{3/2}J^{-1/2}(-\partial_2\log\sqrt{\varepsilon}+\frac{3}{2}c^{-1}\partial_2c+\frac{3}{2}J^{-1}\partial_2J-\frac{1}{2}g^{-1}\partial_2g)(g^{22}I\eta_2+g^{12}I\eta_1)\\
&+2\partial_2\log\sqrt{\varepsilon}g^{-1/2}\varepsilon^{-1/2}c^{3/2}J^{-1/2}(\partial_2g^{22}I\eta_2+g^{22}\partial_2(I\eta_2)+\partial_2g^{12}I\eta_1+g^{12}\partial_2(I\eta_1))\\
&-2c^{-1}\partial_2cg^{-1/2}\varepsilon^{-1/2}c^{3/2}J^{-1/2}\partial_2\log\sqrt{\varepsilon}(g^{22}I\eta_2+g^{12}I\eta_1)+\varepsilon^{-1/2}R.
\end{split}
\end{equation}

\begin{equation}\label{equation25}
\begin{split}
&-g^{-1/2}\partial_2(c^2g^{-1/2}\partial_1E_{02})\\
=&-\varepsilon^{-1/2}c^{3/2}J^{-1/2}g^{-1/2}\partial_2\log\sqrt{\varepsilon}(-\partial_1\log\sqrt{\varepsilon}+\frac{3}{2}c^{-1}\partial_1c+\frac{3}{2}J^{-1}\partial_1J-\frac{1}{2}g^{-1}\partial_1g)(g^{12}I\eta_2+g^{11}I\eta_1)\\
&-\frac{3}{2}\varepsilon^{-1/2}c^{3/2}J^{-1/2}g^{-1/2}c^{-1}\partial_2c\partial_1\log\sqrt{\varepsilon}(g^{12}I\eta_2+g^{11}I\eta_1)\\
&+\frac{1}{2}\varepsilon^{-1/2}c^{3/2}J^{-1/2}g^{-1/2}J^{-1}\partial_2J\partial_1\log\sqrt{\varepsilon}(g^{12}I\eta_2+g^{11}I\eta_1)\\
&-\varepsilon^{-1/2}c^{3/2}J^{-1/2}g^{-1/2}\partial_1\partial_2\log\sqrt{\varepsilon}I\eta^1\\
&-\varepsilon^{-1/2}c^{3/2}J^{-1/2}g^{-1/2}\partial_1\log\sqrt{\varepsilon}(\partial_2g^{12}I\eta_2+g^{12}\partial_2(I\eta_2)+\partial_2g^{11}I\eta_1+g^{11}\partial_2(I\eta_1))\\
&-\varepsilon^{-1/2}c^{3/2}J^{-1/2}g^{-1/2}\partial_2\log\sqrt{\varepsilon}(\partial_1g^{12}I\eta_2+g^{12}\partial_1(I\eta_2)+\partial_1g^{11}I\eta_1+g^{11}\partial_1(I\eta_1))+\varepsilon^{-1/2}R.
\end{split}
\end{equation}

\begin{equation}\label{equation26}
\begin{split}
&\mu^{-1}\partial_2\mu g^{-1/2}c^2g^{-1/2}\partial_1E_{02}\\
=&(-2\partial_2\log\sqrt{\varepsilon}-2c^{-1}\partial_2c)g^{-1/2}c^2g^{-1/2}\partial_1E_{02}\\
=&2\partial_2\log\sqrt{\varepsilon}g^{-1/2}\varepsilon^{-1/2}c^{3/2}J^{-1/2}(-\partial_1\log\sqrt{\varepsilon}+\frac{3}{2}c^{-1}\partial_1c+\frac{3}{2}J^{-1}\partial_1J-\frac{1}{2}g^{-1}\partial_1g)(g^{12}I\eta_2+g^{11}I\eta_1)\\
&+2\partial_2\log\sqrt{\varepsilon}g^{-1/2}\varepsilon^{-1/2}c^{3/2}J^{-1/2}(\partial_1g^{12}I\eta_2+g^{12}\partial_1(I\eta_2)+\partial_1g^{11}I\eta_1+g^{11}\partial_1(I\eta_1))\\
&-2c^{-1}\partial_2cg^{-1/2}\varepsilon^{-1/2}c^{3/2}J^{-1/2}\partial_1\log\sqrt{\varepsilon}(g^{12}I\eta_2+g^{11}I\eta_1)+\varepsilon^{-1/2}R.
\end{split}
\end{equation}

Substituting \eqref{equation23}\eqref{equation24}\eqref{equation25}\eqref{equation26} into \eqref{equation22}, we have

\[
\begin{split}
&\mu g^{-1/2}(H_{12;3}-H_{13;2})\\
=&\mu g^{-1/2}(\partial_3H_{12}-\partial_2H_{13})\\
=&\mu g^{-1/2}\partial_3H_{12}+\varepsilon^{-1/2}c^{3/2}J^{-1/2}g^{-1/2}(-(\partial_2\log\sqrt{\varepsilon})^2-3c^{-1}\partial_2c\partial_2\log\sqrt{\varepsilon}+J^{-1}\partial_2J\partial_2\log\sqrt{\varepsilon})I\eta^2\\
&+\varepsilon^{-1/2}c^{3/2}J^{-1/2}g^{-1/2}\Bigg(-\partial_2\log\sqrt{\varepsilon}\partial_1\log\sqrt{\varepsilon}+\frac{1}{2}\partial_2\log\sqrt{\varepsilon}\partial_1c-\frac{7}{2}\partial_2c\partial_1\log\sqrt{\varepsilon}\\
&\quad\quad\quad\quad\quad\quad+\frac{1}{2}\partial_1J\partial_2\log\sqrt{\varepsilon}+\frac{1}{2}\partial_2J\partial_1\log\sqrt{\varepsilon}\Bigg)I\eta^1\\
&-\varepsilon^{-1/2}c^{3/2}J^{-1/2}g^{-1/2}\partial_2\partial_2\log\sqrt{\varepsilon}I\eta^2-\varepsilon^{-1/2}c^{3/2}J^{-1/2}g^{-1/2}\partial_1\partial_2\log\sqrt{\varepsilon}I\eta^1\\
&-\varepsilon^{-1/2}c^{3/2}J^{-1/2}g^{-1/2}\partial_1\log\sqrt{\varepsilon}(\partial_2g^{12}I\eta_2+g^{12}\partial_2(I\eta_2)+\partial_2g^{11}I\eta_1+g^{11}\partial_2(I\eta_1))\\
&+\varepsilon^{-1/2}c^{3/2}J^{-1/2}g^{-1/2}\partial_2\log\sqrt{\varepsilon}(\partial_1g^{12}I\eta_2+g^{12}\partial_1(I\eta_2)+\partial_1g^{11}I\eta_1+g^{11}\partial_1(I\eta_1))+\varepsilon^{-1/2}R.
\end{split}
\]
Similarly, we also get
\begin{align*}
&\mu g^{-1/2}(H_{11;3}-H_{13;1})\\
=&\mu g^{-1/2}(\partial_3H_{11}-\partial_1H_{13})\\
=&\mu g^{-1/2}\partial_3H_{11}+\varepsilon^{-1/2}c^{3/2}J^{-1/2}g^{-1/2}(-(\partial_1\log\sqrt{\varepsilon})^2-3c^{-1}\partial_1c\partial_1\log\sqrt{\varepsilon}+J^{-1}\partial_1J\partial_1\log\sqrt{\varepsilon})I\eta^1\\
&+\varepsilon^{-1/2}c^{3/2}J^{-1/2}g^{-1/2}\Bigg(-\partial_2\log\sqrt{\varepsilon}\partial_1\log\sqrt{\varepsilon}+\frac{1}{2}\partial_1\log\sqrt{\varepsilon}\partial_2c-\frac{7}{2}\partial_1c\partial_2\log\sqrt{\varepsilon}\\
&\quad\quad\quad\quad\quad\quad+\frac{1}{2}\partial_2J\partial_1\log\sqrt{\varepsilon}+\frac{1}{2}\partial_1J\partial_2\log\sqrt{\varepsilon}\Bigg)I\eta^2\\
&-\varepsilon^{-1/2}c^{3/2}J^{-1/2}g^{-1/2}\partial_1\partial_1\log\sqrt{\varepsilon}I\eta^1+\varepsilon^{-1/2}c^{3/2}J^{-1/2}g^{-1/2}\partial_1\partial_2\log\sqrt{\varepsilon}I\eta^2\\
&-\varepsilon^{-1/2}c^{3/2}J^{-1/2}g^{-1/2}\partial_2\log\sqrt{\varepsilon}(\partial_1g^{12}I\eta_1+g^{12}\partial_1(I\eta_1)+\partial_1g^{22}I\eta_2+g^{22}\partial_1(I\eta_2))\\
&+\varepsilon^{-1/2}c^{3/2}J^{-1/2}g^{-1/2}\partial_1\log\sqrt{\varepsilon}(\partial_2g^{12}I\eta_1+g^{12}\partial_2(I\eta_1)+\partial_2g^{22}I\eta_2+g^{22}\partial_2(I\eta_2))+\varepsilon^{-1/2}R.
\end{align*}
Then
\begin{equation}\label{eqH123combined}
\begin{split}
&\varepsilon^{1/2}c^{-3/2}J^{1/2}g^{1/2}\left(\mu g^{-1/2}(H_{12;3}-H_{13;2})\eta^2+\mu g^{-1/2}(H_{11;3}-H_{13;1})\eta^1\right)\\
=&\varepsilon^{1/2}c^{-3/2}J^{1/2}(\mu\partial_3H_{12}\eta^2+\mu\partial_3H_{11}\eta^1)-\partial_\alpha\log\sqrt{\varepsilon}\partial_\beta\log\sqrt{\varepsilon}\eta^\alpha I\eta^\beta-\partial_\alpha\partial_\beta\log\sqrt{\varepsilon}I\eta^\alpha\eta^\beta\\
&-3c^{-1}\partial_\alpha c\partial_\beta\log\sqrt{\varepsilon}I\eta^\alpha\eta^\beta+J^{-1}\partial_\alpha J\partial_\beta\log\sqrt{\varepsilon}I\eta^\alpha\eta^\beta\\
&-\partial_1\log\sqrt{\varepsilon}(g^{12}\partial_2(I\eta_2)+g^{11}\partial_2(I\eta_1))\eta^2+\partial_2\log\sqrt{\varepsilon}(g^{12}\partial_1(I\eta_2)+g^{11}\partial_1(I\eta_1))\eta^2\\
&-\partial_2\log\sqrt{\varepsilon}(g^{12}\partial_1(I\eta_1)+g^{22}\partial_1(I\eta_2))\eta^1+\partial_1\log\sqrt{\varepsilon}(g^{12}\partial_2(I\eta_1)+g^{22}\partial_2(I\eta_2))\eta^1\\
&-\partial_1\log\sqrt{\varepsilon}(\partial_2g^{12}I\eta_2\eta^2+\partial_2g^{11}I\eta_1\eta^2-\partial_2g^{12}I\eta_1\eta^1-\partial_2g^{22}I\eta_2\eta^1)\\
&-\partial_2\log\sqrt{\varepsilon}(\partial_1g^{12}I\eta_1\eta^1+\partial_1g^{22}I\eta_2\eta^1-\partial_1g^{12}I\eta_2\eta^2-\partial_1g^{11}I\eta_1\eta^1)+\varepsilon^{-1/2}R.\\
\end{split}
\end{equation}
~\\

Now let us put everything together into equation \eqref{mainidentity}.
Using the identities
\[
\partial_kg^{m\ell}=-g^{im}\Gamma^\ell_{ki}-g^{j\ell}\Gamma^m_{kj},
\]
\[
2J^{-1}\partial_\alpha J=g^{\beta\gamma}\partial_\alpha g_{\beta\gamma}=g^{\beta\gamma}(\Gamma_{\alpha\gamma}^\delta g_{\beta\delta}+\Gamma^\delta_{\alpha\beta}g_{\gamma\delta})=2\Gamma^\beta_{\alpha\beta},
\]
we obtain
\begin{equation}\label{identitymiddle1}
\begin{split}
&-2J^{-1}\partial_\alpha J\partial_1\log\sqrt{\varepsilon}\eta^\alpha\eta^1+g^{\alpha 1}J^{-1}\partial_\alpha J\partial_1\log\sqrt{\varepsilon}|\eta|^2\\
&+(\partial_2g^{12}\eta_2\eta^2+\partial_2g^{11}\eta_1\eta^2-\partial_2g^{12}\eta_1\eta^1-\partial_2g^{22}\eta_2\eta^1)\partial_1\log\sqrt{\varepsilon}\\
=&g^{\delta\alpha}\Gamma^1_{\alpha\delta}\partial_1\log\sqrt{\varepsilon}|\eta|^2-2\Gamma_{\alpha\beta}^1\eta^\alpha\eta^\beta\partial_1\log\sqrt{\varepsilon},
\end{split}
\end{equation}
and
\begin{equation}\label{identitymiddle2}
\begin{split}
&-2J^{-1}\partial_\alpha J\partial_2\log\sqrt{\varepsilon}\eta^\alpha\eta^2+g^{\alpha 2}J^{-1}\partial_\alpha J\partial_2\log\sqrt{\varepsilon}|\eta|^2\\
&+(\partial_1g^{12}\eta_1\eta^1+\partial_1g^{22}\eta_1\eta^2-\partial_1g^{12}\eta_2\eta^2-\partial_1g^{11}\eta_2\eta^1)\partial_2\log\sqrt{\varepsilon}\\
=&g^{\delta\alpha}\Gamma^2_{\alpha\delta}|\eta|^2\partial_2\log\sqrt{\varepsilon}-2\Gamma_{\alpha\beta}^2\eta^\alpha\eta^\beta\partial_2\log\sqrt{\varepsilon}.
\end{split}
\end{equation}

Inserting \eqref{eqE113} \eqref{eqE123} \eqref{eqE1312combined} and \eqref{eqH123combined} into the equation \eqref{mainidentity}, using \eqref{identitymiddle1} and \eqref{identitymiddle2}, we obtain

\begin{equation}\label{finalequation}
\begin{split}
\varepsilon^{-1/2}c^{-3/2}J^{1/2}c^{-2}\Bigg(\varepsilon^{-1}\partial_3\varepsilon H_{12}\eta^2+\frac{1}{2}g^{-1}\frac{\partial g}{\partial\tau}H_{12}\eta^2-2\partial_3H_{12}\eta^2\\
+\varepsilon^{-1}\partial_3\varepsilon H_{11}\eta^1+\frac{1}{2}g^{-1}\frac{\partial g}{\partial\tau}H_{11}\eta^1-2\partial_3H_{11}\eta^1\\
+\frac{\sigma}{\varepsilon}H_{12}\eta^2+\frac{\sigma}{\varepsilon}H_{11}\eta^1\\
-g^{1\beta}\frac{\partial g_{\beta1}}{\partial\tau}\frac{1}{\varepsilon}H_{12}\eta^2+g^{1\beta}\frac{\partial g_{\beta2}}{\partial\tau}\frac{1}{\varepsilon}H_{11}\eta^2+g^{2\beta}\frac{\partial g_{\beta1}}{\partial\tau}\frac{1}{\varepsilon}H_{12}\eta^1-g^{2\beta}\frac{\partial g_{\beta2}}{\partial\tau}\frac{1}{\varepsilon}H_{11}\eta^1\Bigg)+IN=0,
\end{split}
\end{equation}
where
\[
\begin{split}
N=&-\partial_{33}^2\log\sqrt{\varepsilon}|\eta|^2+2\partial_\alpha\partial_\beta\log\sqrt{\varepsilon}\eta^\alpha\eta^\beta-g^{\alpha\beta}\partial_\alpha\partial_\beta\log\sqrt{\varepsilon}|\eta|^2-g^{\alpha\beta}\partial_\alpha c\partial_\beta\log\sqrt{\varepsilon}|\eta|^2\\
&-\frac{1}{2}g^{\alpha\beta}\frac{\partial g_{\alpha\beta}}{\partial\tau}\partial_3\log\sqrt{\varepsilon}|\eta|^2+\frac{\partial g_{\alpha\beta}}{\partial\tau}\partial_3\log\sqrt{\varepsilon}\eta^\alpha\eta^\beta\\
&+g^{\alpha\beta}\Gamma^\gamma_{\alpha\beta}\partial_\gamma\log\sqrt{\varepsilon}|\eta|^2-2\Gamma_{\alpha\beta}^\gamma\partial_\gamma\log\sqrt{\varepsilon}\eta^\alpha\eta^\beta\\
&+(\partial_3\log\sqrt{\varepsilon})^2|\eta|^2+g^{\alpha\beta}\partial_\alpha\log\sqrt{\varepsilon}\partial_\beta\log\sqrt{\varepsilon}|\eta|^2+\partial_3\log\sqrt{\varepsilon}\partial_3c|\eta|^2+4c^{-1}\partial_\alpha c\partial_\beta\log\sqrt{\varepsilon}\eta^\alpha\eta^\beta\\
&+\text{terms that do not depend on $\varepsilon$}\\
=&-\nabla_3\nabla_3\log\sqrt{\varepsilon}|\eta|^2-g^{\alpha\beta}\nabla_\alpha\nabla_\beta\log\sqrt{\varepsilon}|\eta|^2\\
&+(\partial_3\log\sqrt{\varepsilon})^2|\eta|^2+g^{\alpha\beta}\nabla_\alpha\log\sqrt{\varepsilon}\nabla_\beta\log\sqrt{\varepsilon}|\eta|^2+4c^{-1}\nabla_\alpha c\nabla_\beta\log\sqrt{\varepsilon}\eta^\alpha\eta^\beta\\
&+\text{terms that do not depend on $\varepsilon$}\\
=&-c^2\Delta\log\sqrt{\varepsilon}+2\nabla_\alpha\nabla_\beta\log\sqrt{\varepsilon}\eta^\alpha\eta^\beta+c^2|\nabla\log\sqrt{\varepsilon}|^2+4c\nabla_\alpha c\nabla_\beta\log\sqrt{\varepsilon}\eta^\alpha\eta^\beta\\
&+\text{terms that do not depend on $\varepsilon$}\\
=:&A^{\alpha\beta}\eta_\alpha\eta_\beta+C.
\end{split}
\]
Here the second order tensor $A$ is of the form
\[
A=-\Delta\log\sqrt{\varepsilon}g_E+2\nabla^2\log\sqrt{\varepsilon}+|\nabla\log\sqrt{\varepsilon}|^2g_E+4c^{-1}\nabla c\otimes\nabla\log\sqrt{\varepsilon},
\]
and $C$ might depend on $\varepsilon\mu,\frac{\sigma}{\varepsilon}$, but not $\varepsilon$.
Recall that $g_E$ is the Euclidean metric.
Here we have used the facts
\[
\nabla_{3}\nabla_{3}\log\sqrt{\varepsilon}=\partial^2_{33}\log\sqrt{\varepsilon}-c^{-1}\partial_3c\partial_3\log\sqrt{\varepsilon}+cg^{\alpha\beta}\partial_\beta c\partial_\alpha\log\sqrt{\varepsilon},
\]
\[
\nabla_{\alpha}\nabla_{\beta}\log\sqrt{\varepsilon}=\partial_\alpha\partial_\beta\log\sqrt{\varepsilon}+\frac{1}{2}c^{-2}\frac{g_{\alpha\beta}}{\partial\tau}\partial_3\log\sqrt{\varepsilon}-\Gamma^\gamma_{\alpha\beta}\partial_\gamma\log\sqrt{\varepsilon}.
\]
Note
\[
\begin{split}
&-\varepsilon^{1/2}c^{-3/2}J^{1/2}\left(g^{1\beta}\frac{\partial g_{\beta1}}{\partial\tau}\frac{1}{\varepsilon}H_{12}\eta^2-g^{1\beta}\frac{\partial g_{\beta2}}{\partial\tau}\frac{1}{\varepsilon}H_{11}\eta^2-g^{2\beta}\frac{\partial g_{\beta1}}{\partial\tau}\frac{1}{\varepsilon}H_{12}\eta^1+g^{2\beta}\frac{\partial g_{\beta2}}{\partial\tau}\frac{1}{\varepsilon}H_{11}\eta^1\right)\\
=&-g_{\alpha\beta}\frac{\partial g_{\alpha\beta}}{\partial\tau}|\eta|^2+\frac{\partial g_{\alpha\beta}}{\partial\tau}\eta^\alpha\eta^\beta.
\end{split}
\]
Denote
\[
H_{1\alpha}=J^{-1/2}c^{3/2}\varepsilon^{1/2}I\xi_\alpha.
\]
Using
\[
\begin{split}
&-2J^{1/2}c^{-3/2}\varepsilon^{-1/2}\partial_3H_{1\alpha}\\
=&J^{-1}\partial_3JI\xi_\alpha-3c^{-1}\partial_3cI\xi_\alpha-\varepsilon^{-1}\partial_3\varepsilon I\xi_\alpha-2\frac{D\xi_\alpha}{d\tau}I-g^{\beta\gamma}\frac{\partial g_{\alpha\beta}}{\partial\tau}I\eta_\gamma+2c^{-1}\partial_3 cI\eta_\alpha+\frac{\sigma}{\varepsilon}I.
\end{split}
\]
Then the equation \eqref{finalequation} can be written as
\[
-2c^{-2}\frac{D\xi_\alpha}{d\tau} \eta^\alpha=-2c^{-2}g^{\alpha\beta} \frac{D\xi_\alpha}{d\tau}\eta_\beta=-2c^{-2}\langle \frac{D\xi}{d\tau},\eta\rangle=-2\frac{D}{d\tau}\left(c^{-2}\langle \xi,\eta\rangle\right)=-2\frac{D}{d\tau}\langle \xi,\eta\rangle_{g}=-N.
\]
Consider the above equation as an ODE along a geodesic $\gamma$ connecting two boundary points $a,b\in\Gamma$.
Then
\[
c^{-2}(\langle\xi,\eta\rangle(b)-\langle\xi,\eta\rangle(a))=\frac{1}{2}\int_\gamma A_{ij}(\gamma(s))\eta^i(s)\eta^j(s) +C\mathrm{d}s.
\]
Here $\eta$ can be any parallel vector field orthogonal to $\gamma$.
 This means that we can recover the (local) transverse ray transform (cf. \cite{Shara,de2021generic,uhlmann2024invertibility1}) of the second order tensor $A$ from the impedance map. We remark here that the same transform has also been derived in \cite{romanov2002inverse} written in terms of $\mu=c^{-2}\varepsilon^{-1}$ instead of $\varepsilon$. One can easily check their equivalence.
By the local invertibility of the transverse ray transform \cite{uhlmann2024invertibility1,uhlmann2024invertibility2}, we can recover the symmetrization of the tensor $A$ in $\kappa^{-1}([0,1])\cap\overline{\Omega}$.

To conclude the above discussions, if $\Lambda_{\varepsilon_1,\mu_1,\sigma_1}=\Lambda_{\varepsilon_2,\mu_2,\sigma_2}$, we have
\[
\begin{split}
&-\Delta\log\sqrt{\varepsilon_1}g_E+2\nabla^2\log\sqrt{\varepsilon_1}+|\nabla\log\sqrt{\varepsilon_1}|^2g_E+4c^{-1}\nabla c\otimes^s\nabla\log\sqrt{\varepsilon_1}\\
=&-\Delta\log\sqrt{\varepsilon_2}g_E+2\nabla^2\log\sqrt{\varepsilon_2}+|\nabla\log\sqrt{\varepsilon_2}|^2g_E+4c^{-1}\nabla c\otimes^s\nabla\log\sqrt{\varepsilon_2}
\end{split}
\]
in $\kappa^{-1}([0,1])\cap\overline{\Omega}$. Here $\otimes^s$ denotes the symmetrized product, i.e.,
\[
u\otimes^s v=\frac{1}{2}(u\otimes v+v\otimes u).
\]
Taking divergence twice of the above equation, we get a fourth order elliptic partial differential equation of the form
\[
\Delta^2(\log\sqrt{\varepsilon_1}-\log\sqrt{\varepsilon_2})+\sum_{|\alpha|\leq 3}F^\alpha(\varepsilon_1,\varepsilon_2)D_\alpha(\log\sqrt{\varepsilon_1}-\log\sqrt{\varepsilon_2})=0\quad\text{in }\kappa^{-1}([0,1])\cap\overline{\Omega}.
\]
Using the boundary determination results (cf. Theorem \ref{boundaryjets})  $\partial_\nu^m(\log\sqrt{\varepsilon_1}-\log\sqrt{\varepsilon_2})=0$, $m=0,1,2,3$ on $\Gamma$, and the \textit{Unique Continuation Principle} (cf. \cite{protter1960unique}) , we end up with
\[
\varepsilon_1=\varepsilon_2\quad \text{in }\kappa^{-1}([0,1])\cap\overline{\Omega}.
\]

\bibliographystyle{abbrv}
\bibliography{biblio}

\end{document}